\documentclass[11pt]{article}
\usepackage{amsmath, amssymb, epsfig, eepic}

\numberwithin{equation}{section}

\newcommand{\eps}{\varepsilon}

\newcommand{\qed}{\hfill $\Box$}
\newcommand{\proof}{\noindent {\bf Proof}. \hspace{2mm}}
\newcommand{\prob}{\mathbb P}
\newcommand{\expec}{\mathbb E}
\newcommand{\smallsup}[1] {{\scriptscriptstyle{({#1}})}}
\newtheorem{theorem}{Theorem}[section]
\newtheorem{lemma}[theorem]{Lemma}
\newtheorem{prop}[theorem]{Proposition}
\newtheorem{cor}[theorem]{Corollary}

\newtheorem{ass}[theorem]{Assumption}
\newtheorem{remark}[theorem]{Remark}

\newcommand{\vep}{\varepsilon}
\newcommand{\eq}{\begin{equation}}
\newcommand{\en}{\end{equation}}
\def\eqalign#1\enalign{
    \begin{align}#1\end{align}
    }
\newcommand{\sss}   { \scriptscriptstyle }
\newcommand{\nn}   { \nonumber}
\newcommand{\Q}{\mathbb{Q}}

\newcommand{\cZ}{{\cal Z}}
\newcommand{\CN}{{\mathbf C}_{\sN}}

\def\NN{{\mathbb N}}

\def\indic{{I}}

\newcommand*{\sumu}{\displaystyle\sum}
\newcommand*{\produ}{\displaystyle\prod}
\def\AeN{ {\cal A}_{\varepsilon\kern-.05em,\kern .1em N} }
\def\BeN{ B_{\varepsilon\kern-.05em,\kern .1em N}}
\def\CeN{ {\cal C}_{\varepsilon\kern-.05em,\kern .1em N} }
\def\DeN{ {\cal D}_{\varepsilon\kern-.05em,\kern .1em N} }
\def\EeN{ {\cal E}_{\varepsilon\kern-.05em,\kern .1em N} }

\def\prob{{\mathbb P}}
\def\expec{{\mathbb E}}
\newcommand{\sN}{{\sss N}}
\newcommand{\sZ}{{\sss Z}}
\newcommand{\sF}{{\sss F}}
\newcommand{\Ccal}{{\cal C}}
\newcommand{\whp}{{\bf whp} }
\newcommand{\whps}{{\bf whp}}

\newcommand*{\picdirrectory}{}
\newcommand*{\fig}[3]{
        \begin{figure}[!h!t]
    \begin{center}
        \includegraphics{\picdirrectory#1}
    \end{center}
        \caption{\footnotesize #2}        \label{#3}
        \end{figure}}

\usepackage[usenames]{color}

\newcommand{\Core}{{\rm Core}_{\sN}}
\newcommand{\Ncal}[1]{{\cal N}^{\sss(#1)}}

\newcommand{\N}{{\mathbb N}}

\title{Random graphs with arbitrary i.i.d.\ degrees}
\author{Remco van der Hofstad\footnote{Department of Mathematics and
Computer Science, Eindhoven University of Technology, P.O.\ Box
513, 5600 MB Eindhoven, The Netherlands. E-mail: {\tt
rhofstad@win.tue.nl}}\\
Gerard Hooghiemstra\footnote{Delft University of Technology,
Electrical Engineering, Mathematics and Computer Science, P.O. Box
5031, 2600 GA Delft, The Netherlands. E-mail: {\tt
G.Hooghiemstra@ewi.tudelft.nl}}$\,\,$ and Dmitri
Znamenski\footnote{EURANDOM, P.O.\ Box 513, 5600 MB Eindhoven, The
Netherlands. E-mail: {\tt znamenski@eurandom.nl}} }

\begin{document}
\maketitle

\begin{abstract}
In this paper we derive results concerning the connected components
and the diameter of random graphs with an arbitrary i.i.d.\ degree
sequence. We study these properties primarily, but not exclusively, when the tail of the degree
distribution is regularly varying with exponent $1-\tau$. There are
three distinct cases: (i) $\tau>3$, where the degrees have finite
variance, (ii) $\tau\in (2,3)$, where the degrees have infinite
variance, but finite mean, and (iii) $\tau\in (1,2)$, where the
degrees have infinite mean. These random graphs can serve as models
for complex networks where degree power laws are observed.

Our results are twofold. First, we give a
criterion when there exists a unique largest connected component of
size proportional to the size of the graph, and study sizes of the
other connected components.
Secondly, we establish a phase transition for the diameter when $\tau \in (2,3)$.
Indeed, we show that for $\tau>2$ and when nodes with degree
$2$ are present with positive probability, the diameter of the random graph is, with high probability,
bounded below by a constant times the logarithm of the size of the graph.
On the other hand, assuming that
all degrees are at least 3 or more, we show that, for $\tau \in (2,3)$,
the diameter of the graph is with high probability bounded from above by
a constant times the $\log \log$ of the size of the graph.

\end{abstract}


\section{Introduction}
Random graph models for complex networks have received
a tremendous amount of attention in the past decade.
Measurements have shown that many real networks share two
properties. The first fundamental network property is the fact
that typical distances between nodes are small. This is called
the `small world' phenomenon (see \cite{Watt99}).
For example, in the Internet, IP-packets cannot use more
than a threshold of physical links, and if the distances in
terms of the physical links would be large, e-mail service would simply break down.
Thus, the graph of the Internet has evolved in such a way that
typical distances are relatively small, even though the Internet
is rather large. The second and maybe more surprising property
of many networks is that the number of nodes with degree $k$ falls
off as an inverse power of $k$. This is called a `power law degree
sequence', and resulting graphs often go under the name
`scale-free graphs', which refers to the fact that the asymptotics of the
degree sequence is independent of the size of the graph (see \cite{FFF99}).
We refer to \cite{AB02, Newm03, Stro01} and the references therein for an
introduction to complex networks and many examples where the
above two properties hold.

The observation that many real networks
have the above two properties has incited a burst of activity in network
modeling using random graphs. These models can be divided into two distinct types:
`static' models, where we model a graph of a given size as a time
snap of a real network, and `dynamical' models, where we model the growth
of the network. Static models aim to describe real networks
and their topology at a given time instant.
Dynamical models aim to {\it explain} how the networks came to be as they are. Such explanations often
focus on the growth of the network as a way to explain the
power law degree sequences by means of `preferential attachment'
growth rules, where added nodes and edges are more likely to be
attached to nodes that already have large degrees. See \cite{Bara02}
for a popular account of preferential attachment.

The random graph where the degrees are i.i.d.\
is sometimes called the {\it configuration model} (see \cite{Newm03}).
In this paper, we study properties
of the connected components in the random graph with i.i.d.\
degrees, and prove results
concerning the scaling of the largest and second largest
connected components, as well as the diameter.

The remainder of this introduction is organized as follows.
In Section \ref{sec-mod} we start by introducing the configuration model, and
in Section \ref{sec-connres}, we discuss the new results concerning
component sizes and diameter of this graph.
We describe related work
and open questions in Section \ref{sec-OP}. We complete the
introduction with the organization of the paper in Section \ref{sec-org}.


\subsection{The configuration model}
\label{sec-mod}
Fix an integer $N$. Consider an i.i.d.\ sequence of random variables
$D_1,D_2,\ldots,D_{\sN}$. We will construct an undirected graph
with $N$ nodes where node $j$ has degree $D_j$. We will assume that
$L_{\sN}=\sum_{j=1}^N D_j$ is even. If $L_{\sN}$ is odd, then we
add a stub to the $N^{\rm th}$ node, so that $D_{\sN}$ is
increased by 1. This single change will make hardly any difference
in what follows, and we will ignore this effect. We will later
specify the distribution of $D_1$.

To construct the graph, we have $N$
separate nodes and incident to node $j$, we have $D_j$ stubs or half-edges. All stubs
need to be connected to build the graph. The stubs are
numbered in a given order from $1$ to $L_{\sN}$. We start by
connecting at random the first stub with one of the $L_{\sN}-1$
remaining stubs. Once paired, two stubs form a single edge of the
graph. Hence, a stub can be seen as the left or the right half of
an edge. We continue the procedure of randomly choosing and
pairing the stubs until all stubs are connected.
Unfortunately, nodes having self-loops may occur. However,
self-loops are scarce when $N \to \infty$, as shown in \cite{BDM-L05}.

The above model is a variant of the {\it configuration model},
which, given a degree sequence, is the random graph with that
given degree sequence. The degree sequence of a graph
is the vector of which the $k^{\rm th}$ coordinate equals the
fraction of nodes with degree $k$. In our model, by the law of large numbers,
the degree sequence is close to the distribution of the nodal degree
$D$ of which $D_1, \ldots, D_{\sN}$ are i.i.d.\ copies.

The probability mass function and the distribution function of the nodal degree
law are denoted by
    \begin{equation}
    \label{kansen}
    \prob(D_1=k)=f_k,\quad k=1,2,\ldots, \quad \mbox{and} \quad
    F(x)=\sum_{k=1}^{\lfloor x \rfloor} f_k,
    \end{equation}
where $\lfloor x \rfloor$ is the largest integer smaller than or
equal to $x$. We pay special attention to distributions of the form
    \begin{equation}
    \label{distribution}
    1-F(x)=x^{1-\tau}L(x),
    \end{equation}
where $\tau >1$ and $L$ is slowly varying at infinity. This means that
the random variables $D_j$ obey a power law, and the factor $L$ is meant to
generalize the model. For one of our main results (Theorem \ref{thm-clusterP}
below) we assume the following more specific conditions,
splitting between the cases $\tau\in (1,2), \tau\in (2,3)$ and $\tau>3$:

\begin{ass}
\label{ass}
\begin{enumerate}
\item[{\rm (i)}] For $\tau\in (1,2)$, we assume (\ref{distribution}).
\item[{\rm (ii)}] For $\tau \in (2,3)$, we assume that
there exists $\gamma\in [0,1)$ and $C>0$ such that
    \eq
    \label{Fcond}
    x^{1-\tau-C(\log{x})^{\gamma-1}}\leq 1-F(x)\leq
    x^{1-\tau+C(\log{x})^{\gamma-1}},\qquad \mbox{for large $x$}.
    \en
\item[{\rm (iii)}] For $\tau>3$, we assume that
there exists a constant $c>0$ such that
    \begin{equation}
    \label{tau>3ass}
    1-F(x)\leq c x^{1-\tau},\qquad \mbox{for all $x\geq 1$},
    \end{equation}
and that $\nu>1$, where $\nu$ is given by
    \eq
    \label{nu}
    \nu=\frac{\expec[D_1(D_1-1)]}{\expec[D_1]}.
    \en\end{enumerate}
\end{ass}
Distributions satisfying (\ref{tau>3ass}) include
distributions which have a lighter tail than a power law,
and (\ref{tau>3ass}) is only slightly stronger than
assuming finite variance.
The condition in (\ref{Fcond}) is slightly stronger than
(\ref{distribution}).

\subsection{Connected components and diameter of the random graph}
\label{sec-connres}
In this paper, we prove results concerning the sizes of the
connected components in the random graph, and give bounds on the
diameter.

For these results, we need some additional notation.
For $\tau>2$, we introduce a
{\it delayed} branching process $\{ {\cal Z}_n\}_{n\geq 1}$, where in the
first generation the offspring distribution is chosen according to
(\ref{kansen}) and in the second and further generations the offspring is
chosen in accordance to $g$ given by
    \begin{equation}
    \label{outgoing degree}
        g_k=\frac{(k+1) f_{k+1}}{\mu},\quad k=0,1,\ldots,
    \qquad\text{where}\qquad
    \mu=\expec[D_1].
    \end{equation}
%

In the statements below, we write $G$ for the random graph
with degree distribution given by (\ref{kansen}), and we denote
for $\tau>2$ the survival probability of the delayed branching
process $\{{\cal Z}_n\}$ described above by $q$.
When $1<\tau<2$, for which $\mu=\expec[D_1]=\infty$, we define $q=1$.
We define, for $\delta>0$,
    \eq
    \label{gammadef}
    \gamma^*_1=\frac{1+\delta}{\log{\mu}-\log{2}}\quad (\tau>2), \qquad
    \gamma^*_2=\frac{\tau-1}{2-\tau}(1+\delta)\quad (\tau\in (1,2)).
    \en
In the sequel we use the abbreviation
{\bf whp} to denote that a statement holds with probability $1-o(1)$
as $N\rightarrow \infty$.

\begin{theorem}[The giant component]
    \label{thm-clusterP}
    Fix $\delta>0$.
    When Assumption \ref{ass} holds and $q\in (0,1]$, then,
    {\bf whp}, the largest connected component in $G$ has
    $qN(1+o(1))$ nodes, and all other connected components have at most
    $\gamma^*_2$ nodes when $\tau\in (1,2)$, and at most $\gamma^*_1 \log{N}$
    nodes when $\tau>2$ and $\mu>2$.
\end{theorem}

Theorem \ref{thm-clusterP} is similar in spirit to the main results in
\cite{MR95, MR98}, where the connected components in the configuration
model were studied for {\it fixed} degrees, rather than i.i.d.\ degrees.
In \cite{MR95, MR98}, however, restrictions were posed on the maximal
degree. Indeed, for the asymptotics of the largest connected component,
it was assumed that the maximal degree is bounded by $N^{\frac 14 -\vep}$,
for some $\vep>0$. Since, for i.i.d.\ degrees, the maximal degree is
of the order $N^{\frac{1}{\tau-1}+o(1)}$, where $\tau>1$ is the degree exponent,\
this restricts to $\tau>5$. Theorem \ref{thm-clusterP} allows for any $\tau>1$,
at the expense of the assumption that $\mu>2$ and the fact that the degrees
are i.i.d. The latter restriction comes from the fact that,
in the proofs, we make essential use of the results in \cite{EHHZ04, HHV05, HHZ04a}.
However, a close inspection of the proofs in \cite{EHHZ04, HHV05, HHZ04a} shows that
independence is not exactly needed. This is explained in more detail for the
case that $\tau\in(2,3)$ in \cite{HHZ04a}. The restriction $\mu>2$ is somewhat unusual,
and is not present in \cite{MR95, MR98}. However, for most real networks,
this condition {\it is} satisfied (see e.g., \cite{BA99, Newm03}).

In \cite{HHV05}, a similar result as Theorem \ref{thm-clusterP} was proved for the case when $\tau>3$, using
the results of \cite{MR95, MR98}, without the assumption that $\mu>2$.
In this case, the main restriction is that $\nu>1$. This result is proved
by suitably adapting the graph by erasing some edges from the nodes
with degree larger than $N^{\frac14-\eps}$.
The result in \cite{MR95, MR98} can be restated as saying that
the largest component is $qN(1+o(1))$ when $\nu>1$ and is $o(N)$ when
$\nu<1$. Our result applies in certain cases where the results of
\cite{MR95, MR98} do not apply (such as the cases when
$\tau\in (2,3)$ and $\tau\in(1,2)$), and our proof is relatively
simple and yields rather explicit bounds.

The proof of Theorem \ref{thm-clusterP} is organized as follows. In
\cite{HHV05} and \cite{HHZ04a}, respectively,
it was shown that for $\tau >2$ the probability that two nodes are
connected is asymptotically equal to $q^2$, where
$q$ arises as the survival probability of the branching process
approximation of the shortest-path graph from a given node. For
$\tau \in (1,2)$ this branching process is not defined and we use the convention $q=1$, because for
$\tau \in (1,2)$  it was shown in \cite{EHHZ04} that the probability that two arbitrary nodes are connected equals 1 with high probability.
These results suggest that there exists a largest connected component
of size roughly equal to $qN$.
The proofs in \cite{HHV05,HHZ04a} rely on branching process comparisons
of the number of nodes that can be reached within $k$ steps.
The proof in \cite{EHHZ04} relies mainly on extreme value theory.
The main ingredient in the proof of Theorem \ref{thm-clusterP}
is that we show that, when $\mu>2$, any connected component is
either very large, or bounded above by $\gamma^*_2$ when
$\tau\in (1,2)$, and by $\gamma^*_1 \log{N}$ when $\tau>2$
(see Proposition \ref{thm-cs} below). The proof of these facts
again relies on branching process comparisons, using the detailed
estimates obtained in \cite{HHV05, HHZ04a}. Since any two nodes are
connected to each other with positive probability, there must be
at least one such large connected component. The proof is completed
by showing that this largest connected component of size proportional to
the size of the graph is unique, and that its size is close to
$qN$.

While Theorem \ref{thm-clusterP} provides good upper bounds on the second
largest component and detailed asymptotics on the largest component,
it leaves a number of questions open. For example, how large is
the second largest component, and, when $q=1$, is the graph connected?
We next investigate these questions.

The following theorem says that $\gamma^*_1$ and $\gamma^*_2$ defined
in~(\ref{gammadef}) provide quite sharp estimates for the size of
the connected components that are not the
largest in the random graph. Define for $f_1>0$,
    \eq
    \label{gamma**}
    \gamma^{**}_1=\frac{1-\delta}{\log\mu-\log f_1}\qquad (\tau>2),
    \qquad
    \gamma^{**}_2=\frac{\tau-1}{2-\tau}(1-\delta) \qquad (\tau\in (1,2)).
    \en
\begin{theorem}[Sizes of non-giant components]
\label{dnz9Th2}
\begin{itemize}
\item [{\rm (i)}]
Let $\tau\in(1,2)$ and $f_1>0$. Then, for any $\delta>0$
and $k\le \gamma^{**}_2$, and such that $f_k>0$,
\whp the random graph contains a connected component with $k+1$
nodes.

\item[{\rm (ii)}]
Let $\tau>2$ and $\mu>f_1>0$, and assume that
    \begin{equation}
    \label{dnz9_10}
    f_k=L_f(k)k^{-\tau},\quad k\to \infty,
    \end{equation}
where $L_f(\cdot)$ is a slowly varying function.
Then, for any $\delta>0$, and $k=k_{\sN}\le\gamma^{**}_1\log{N}$
and such that $f_k>0$, \whp the random graph contains a
connected component with $k+1$ nodes.
\end{itemize}
\end{theorem}

We present some further results in the more special case when $q=1$. In this case,
either $\mu=\infty$ or $D\geq 2$ a.s. Then, from Theorem \ref{thm-clusterP}, we have
that there exists a unique connected component of size $N-o(N)$ and all other connected
components are much smaller. For this case we investigate when the random graph is
\whp connected. Let $\CN$ denote the number of nodes in the complement of the
largest connected component of the random graph.

\begin{theorem}[Size of complement of giant component]
\label{dnzTh1intro}
\begin{itemize}
\item[{\rm (i)}] Let $\prob(D\ge 2)=1$ and $2<\mu<\infty$. Then,
there exists $a<1$ and $b>0$ such that for each $1\le k\le N$, as $N\to \infty$,
    \eq
    \label{119}
    \prob(\CN\geq k)\leq b a^{k}.
    \en

\item[{\rm (ii)}] If  $\prob(D\ge 3)=1$, then
    \eq
    \label{conninfmean}
    \lim_{N\rightarrow \infty} \prob(\CN=0)=1.
    \en
Consequently, in the latter case, the random graph is connected
{\bf whp}.

\item[{\rm (iii)}] The conclusion (\ref{conninfmean}) also holds when $\mu=\infty$ and
$\prob(D\ge 2)=1$ instead of $\prob(D\ge 3)=1$.

\item[{\rm (iv)}] The conclusion (\ref{conninfmean}) also holds when
$L_{\sN}/N^2\to \infty$ in probability
without any further restrictions on the degree distribution.
\end{itemize}
\end{theorem}

Clearly, the restriction $\prob(D\ge 2)=1$ is necessary to obtain a connected graph, but it is not sufficient,
as we will see below. Equation \eqref{119}
establishes that the complement of the connected component has exponential tails.
When $f_2=\prob(D=2)>0$, then it is not hard to see that the expected number of
pairs of nodes with degree equal to two that are connected to each other, for $\mu < \infty$,
is asymptotically equal to
    \eq
    \label{2cycles}
    \frac{(f_2N)^2}{2}\frac{2}{(\mu N)^2}=\left(\frac{f_2}{\mu}\right)^2>0.
    \en
Indeed, $(Nf_2)^2/2$ is roughly equal to the
number of pairs of nodes with degree 2, and $2(L_{\sN})^{-2}$ is roughly
equal to the probability that the two stubs between these two nodes are
connected to each other. The above mean is strictly positive, which suggests
that the number of
pairs of nodes with degree equal to two that are connected to each other
is with strictly positive probability positive. We believe that
the proof in \cite[Section 2.4]{Boll01} can be followed to show that
this number is close to a Poisson distribution with parameter $(f_2/\mu)^2$.
Similar computations can be performed for the number of cycles of length
3 or larger consisting of nodes with degree precisely equal to 2. Thus, for
$f_2>0$, \eqref{119} seems the best possible result.
We show in Theorem \ref{dnzTh1intro}(ii) that,
when $f_1+f_2=0$, the graph is connected {\bf whp}. The same result holds
(see Theorem \ref{dnzTh1intro} (iii) and (iv)) when $f_1=0$ and
$\mu=\infty$, or when $L_{\sN}/N^2\to \infty$ in probability.

Finally, we give, in Theorem \ref{thm-diameterfirst} and
\ref{thm-diametersec} below, bounds on the diameter of the graph,
which we define as the largest distance between
any two nodes that are connected:

\begin{theorem}[Lower bound on diameter]
\label{thm-diameterfirst}
For $\tau>2$, assuming that $f_1+f_2>0$ and $f_1<1,$ there
exists a positive constant $\alpha$ such that \whp
the diameter of $G$ is bounded below by $\alpha \log{N}$,
as $N\rightarrow \infty$.
\end{theorem}

The result in Theorem \ref{thm-diameterfirst} is most interesting in the case
when $\tau\in (2,3)$. Indeed, by \cite[Theorem 1.2]{HHZ04a}, the typical distance
for $\tau\in (2,3)$ is proportional to $\log{\log{N}}$, whereas we show here that the diameter is
bounded below by a constant times $\log{N}$ when $f_1+f_2>0$ and $f_1<1$. Therefore, we see that the
average distance and the diameter are of a different order of magnitude,
which is rather interesting. The pairs of nodes where
the distance is of the order $\log{N}$ are thus scarce. The proof of
Theorem \ref{thm-diameterfirst} reveals that these pairs are along
long lines of vertices with degree 2 that are connected to each other.

We end with a theorem stating that when $\tau\in (2,3)$, the above assumption that
$f_1+f_2>0$ is necessary and sufficient for $\log{N}$ lower bounds on the diameter.
We assume that there exists a $\tau\in (2,3)$ such that, for some $c>0$ and all $ x\geq 1$,
    \eq
    \label{infvarass}
    1-F(x) \geq c x^{1-\tau}.
    \en
Observe that (\ref{infvarass}) is strictly weaker than \eqref{Fcond}.
Then the main result is as follows:

\begin{theorem}[Upper bound on diameter]
\label{thm-diametersec}
Assume that $f_1+f_2=0$ and that
\eqref{infvarass} holds. Then, there
exists a positive constant $C_{\sF}$ such that \whp
the diameter of $G$ is bounded above by $C_{\sF} \log\log{N}$,
as $N\rightarrow \infty$.
\end{theorem}

In the course of the proof of Theorem \ref{thm-diametersec}, we will
establish an explicit expression for $C_{\sF}$ in terms of $F$.

We remark that Theorems \ref{dnz9Th2}--\ref{thm-diametersec} do not rely on
Assumption \ref{ass}, while Theorem \ref{thm-clusterP} {\it does}.
The reason for this is that the proof of Theorem \ref{thm-clusterP}
relies on the results proved in \cite{EHHZ04, HHV05, HHZ04a}, while the proofs
of Theorems \ref{dnz9Th2}--\ref{thm-diametersec} are completely self-contained.

\subsection{Related results and open problems in static models}
\label{sec-OP}
As mentioned in the introduction the results in this paper are partly based
on a coupling between the configuration model and branching processes, presented in two
previous publications \cite{HHV05} and \cite{HHZ04a}.
For later reference, we will summarize the graph distance results obtained in these papers
and in \cite{EHHZ04}, for the case $\tau\in(1,2)$.

The graph distance $H_{\sN}$
between the nodes $1$ and $2$ is defined as the minimum number of edges that
form a path from $1$ to $2$. By convention, the distance equals $\infty$ if $1$
and $2$ are not connected. Observe that the distance
between two randomly chosen nodes is equal in distribution to
$H_{\sN}$, because the nodes are exchangeable.
The main result in  \cite{EHHZ04} is that for $\tau \in (1,2)$  and in the limit for
$N$ tending to infinity, the distribution of the graph distance is concentrated
on the points $2$ and $3$, i.e.,
when Assumption \ref{ass} holds, then,
\begin{equation}
\label{limit law22}
\lim_{N\to \infty}
\prob(H_{\sN} =2)=1-\prob(H_{\sN} =3)=p,
\end{equation}
where $p=p_{F}\in (0,1)$.
For $\tau \in (2,3)$ we showed in \cite{HHZ04a}, that when
Assumption \ref{ass} holds, the fluctuations
of $H_{\sN}$ around
\eq
\label{centering23}
2\frac{\log\log N}{|\log (\tau -2)|}
\en
are $O_p(1)$ as $N\to \infty$; and finally, we showed in \cite{HHV05} that the same result holds
for $\tau>3$, with the centering in (\ref{centering23}) replaced by
\eq
\label{centering>3}
\log_{\nu} N.
\en

The model studied in this paper with $\tau\in (2,3)$ is also studied in
\cite{RN04}, where it is proved that {\bf whp} the graph distance
$H_{\sN}$ is less than
$2\frac{\log \log N}{|\log(\tau-2)|} +2\kappa(N)$, where
    \eq
    \kappa(N)=
    \Big\lceil \exp{\Big(\frac{2}{3-\tau} \ell(N)\Big)}\Big\rceil
    \qquad\text{with}\qquad
    \lim_{N\rightarrow \infty} \frac{\ell(N)}{\log\log\log\log N}
    =\infty.
    \en
At approximately the same moment the
$\log\log N$-scaling result appeared in the physics literature \cite{CH03},
where it was derived in a non-rigorous way.
The distance results in \cite{HHV05} was generalized to a much larger class of
random graphs in \cite{EHH06}.

There is substantial work on random graphs that are, although different
from ours, still similar in spirit.
In \cite{ACL01a,DGM02,DMS03,MR95,MR98}, random graphs
were considered with a degree sequence that is {\it precisely}
equal to a power law, meaning that the number of nodes with degree
$k$ is precisely proportional to $k^{-\tau}$.
A second related model can be found in \cite{CL02a, CL02b}, where edges
between nodes $i$ and $j$ are present with probability
equal to $w_iw_j/\sum_l w_l$ for some `expected degree vector'
$w=(w_1, \ldots, w_{\sN})$.
In \cite{CL04}, these authors study  a so-called {\it hybrid} model.

Arratia and Liggett \cite{AL04} study whether {\it simple}
graphs exist with an i.i.d.\ degree distribution, i.e., graphs
without self-loops and multiple edges. It is not hard to see that
when $\tau<2$ this happens with probability 0 (since the largest degree
is larger than $N$). When $\tau>2$, however, this probability is asymptotic
to the probability that the sum of $N$ i.i.d.\ random variables is even, which
is close to $1/2$. When $\tau=2$, the probability can converge to any
element of $[0,\frac 12]$, depending on the slowly varying function
in (\ref{distribution}). A similar problem is addressed in \cite{BDM-L05},
where various ways how self-loops and multiple edges can be
avoided are discussed. Among others, in \cite{BDM-L05}, it is proved that when
the degrees are i.i.d.\ and all self-loops and multiple edges are removed,
then the power law degree sequence remains valid.

There are many open questions remaining in the configuration model.
For instance, in \cite{HHV05}, we have shown that for $\tau>3$,
the largest connected component has size $qN$,
where $q$ is the survival probability of the delayed branching process.
All other connected components have size at most $\gamma \log{N}$,
for some $\gamma>0$. For $\tau\in (2,3)$, such a result is given in this paper under the extra
assumption that $\mu>2$. It would be of interest to investigate
whether the same result holds for $\tau\leq 3$ and general $\mu$ when $q>0$.


A second quantity of interest is the {\it diameter}
of the graph which is important
in many applications. For instance, in the Internet, a message is
killed when the number of hops exceeds a finite threshold. Thus, it
would be interesting to investigate how the diameter grows with the size
of the graph. The result in Theorem \ref{thm-diameterfirst} is a lower bound
in the case when $f_2>0$, whereas Theorem \ref{thm-diametersec} gives an upper bound
for $\tau \in (2,3)$, when $f_1+f_2=0$; however a better understanding of the diameter
is necessary.


An important property of the topology of a graph is its {\it clustering},
which basically describes how likely two nodes that have an edge to a common
node are to be connected by an edge. In general, in random graphs, this
clustering is much smaller than the clustering in real networks. It would
be of interest to investigate graphs with a higher clustering in more detail.
The hybrid graphs in \cite{CL04} are an important step in that direction.



\subsection{Organization of the paper}
\label{sec-org}
The paper is organized as follows.
In Section \ref{sec-conn}, we prove Theorem \ref{dnzTh1intro},
and in Section \ref{sec-conncom}, we prove Theorem \ref{thm-clusterP}.
Finally, in Section \ref{sec-lb}, we prove the lower bounds on
the second largest connected component in Theorem \ref{dnz9Th2} and
on the diameter in Theorem \ref{thm-diameterfirst} and \ref{thm-diametersec} .


\section{Connectivity properties}
\label{sec-conn}
In this section  we prove connectivity properties of the random
graph defined in Section \ref{sec-mod}. In particular we will
prove Theorem \ref{dnzTh1intro}, which states among other things that $\prob(\CN \ge k)$,
where $\CN$ denotes the number of nodes in the complement to the
largest connected component of the random graph, is exponentially bounded
as $N\to \infty$, when $\prob(D\ge 2)=1$ and $2<\mu<\infty$.
Throughout the paper, we write $I[E]$ for the indicator of the event $E$.

We start by stating a lemma which bounds the conditional probability
$\prob_{\sN}(\CN\geq s)$, where $\prob_{\sN}$ denotes
the probability given the degrees $D_1, \ldots, D_{\sN}$.

\begin{prop}
\label{dnzLeTh1} Let $r\in \{1,2\}$, and assume that $\prob(D_1\ge r)=1$. Then,
for any $1\le s\le N/3$,
\begin{equation}
\label{dnz25}
\prob_{\sN}\left(\CN\ge s\right)\le
2\sumu\limits_{j=s}^{N-s}
\left(\frac{2N^{\frac2r}}{L_{\sN}}\right)^{\lceil jr/2\rceil}, \quad a.s.
\end{equation}
\end{prop}

\noindent
We first show that Theorem \ref{dnzTh1intro}(i), (iii) and (iv)
are an immediate consequence of Proposition \ref{dnzLeTh1}.
Theorem \ref{dnzTh1intro}(ii) is proved in Section \ref{sec-loglog} below.

\noindent {\bf Proof of Theorem \ref{dnzTh1intro}(i), (iii) and (iv).}
We start with case (i) where $2<\mu<\infty$, and $\prob(D\geq 2)=1$. We denote
by   $\mu_{\sN}=L_{\sN}/N$. Taking
expectations on both sides of (\ref{dnz25}), yields for $r=2$, and
with $1\le s\le N/3$,
    \begin{align}
    \label{dnz28}
    \prob(\CN\geq s)&\leq
    \expec\left(2\sumu\limits_{j=s}^{N-s}
    \left(\frac{2}{\mu_{\sN}}\right)^j\indic[\mu_{\sN}\ge1+\mu/2]\right)+
    \prob\left(\mu_{\sN}<1+\mu/2\right)\\
    &\le
    \frac{2\left(\frac{2}{1+\mu/2}\right)^{s}}{1-\frac{2}{1+\mu/2}}+
    \prob\left(\mu_{\sN}<1+\mu/2\right)
    \le\frac{2(2+\mu)}{\mu-2}\left(\frac{4}{2+\mu}\right)^s +e^{-I N},
    \nonumber
    \end{align}
where $I$ is the exponential rate of the event $\mu_{\sN}=L_{\sN}/N<(1+\mu/2)$, which is
strictly positive since $\{1+\mu/2<\mu\}$.
Indeed, the final inequality of (\ref{dnz28}) holds for all $N\ge 1$, because
of Chernov's bound and the fact that $\mu>2$,
and that for $t>0$ the Laplace transform $\expec[\exp\{-tD_1\}]$ exists.
For $2<\mu<\infty$, and $\prob(D\geq 2)=1$, this shows that (\ref{119}) holds for all
$N\ge 1$ and $1\le s\le N/3$, by taking $a=\max\{e^{-3I},\frac{4}{2+\mu }\}$,
and $b=1+\frac{2(2+\mu)}{\mu-2}$.
The statement for $N/3\le s \le N$ follows from:
\eq
\label{n3totn}
\prob(\CN\geq s)\leq
\prob(\CN\geq N/3)\leq
ba^{N/3}
=
b(a^{1/3})^N\le b(a^{1/3})^s.
\en
This proves (i).

Consider next the case (iii), where $\mu=\infty$ and $\prob(D\geq 2)=1$. Then,
for any $\varepsilon>0$ and large enough $N$,
    $$
    \prob\left(\frac{2N}{L_{\sN}}\le\varepsilon\right)\ge1-\varepsilon.
    $$
Hence, the probability that the random graph is
disconnected, or $\prob(\CN\ge1)$, is due to~(\ref{dnz25}) at
most
    \eq
    \expec\left(2\sumu\limits_{j=1}^{N-1}
    \left(\frac{2N}{L_{\sN}}\right)^j\indic[2N\le\varepsilon L_{\sN}]\right)+
    \prob\left(\frac{2N}{L_{\sN}}>\varepsilon\right)
    \le2\sumu\limits_{j=1}^{N-1}
    \varepsilon^j+\varepsilon<\frac{2\varepsilon}{1-\varepsilon}+\varepsilon.
    \en
Since $\varepsilon>0$ can be chosen arbitrarily small the
probability that the graph is disconnected tends to zero,
as $N$ tends to infinity.

We complete the proof with the case (iv), where $L_{\sN}/N^2\to\infty$, in probability, which for example is the case when
$\tau\in (1,\frac 32)$.
In this case, we take $r=s=1$, and we use the assumption
that $N^2/{L_{\sN}}\to0$ in probability as $N\to\infty$, to see that
    \eq
        \prob(\CN\ge 1)\leq
        2\sumu\limits_{j=1}^{N-1} (2\vep)^{\lceil j/2\rceil}+\prob(N^2/{L_{\sN}}>\vep)
    \leq 9\vep,
    \en
for each $\vep>0$.
\qed
\vskip0.5cm

\noindent
{\bf Proof of Proposition~\ref{dnzLeTh1}.}
For any $s\leq N/3$, we estimate the probability that $\CN\ge s$.
If $\CN\ge s$, then there exist two disjoint sets of nodes, one
with $s\le j\le N-s$ nodes and another with $N-j$,
such that all stubs of the first set pair within the first set
and all stubs of the second set pair within the second.
To see this, we note that when the largest component has size at least $s$,
then the statement is correct, and the two disjoint sets are the nodes in the
largest component and the ones outside of the largest component.
Thus we are left to prove the statement when the largest component
has size at most $s-1$. In this case we order the connected components by size
(and when there are multiple components of the same size, we do so in an arbitrary way).
Then we start with the largest component and we successively add to it the largest
component that is still available, until the total size of these connected components
is larger than $s$. Since the largest component has size at most $s-1$, the total
size we end up with is in between $s$ and $2s$. Since $s\leq N/3$
we obtain $2s \leq N-s$, and we arrive at the claim that the set consists of a number of nodes
which is in between $s$ and $N-s$. Put the chosen connected components
into the first set  and all remaining nodes in the second. By construction there are no
edges between these two sets of nodes.

Our plan is to show that the probability that the first set of size $j$,
where $s\le j \le N-s$, pairs within its own group and the second set
of size $N-j$ also pairs within its own group is bounded by the right side
of (\ref{dnz25}). We use Boole's inequality to bound the probability
of the union of all possible choices for the first set.

To this end, let $i_1,\dots,i_j$ be the nodes and $A_j=\sum_{l=1}^j D_{i_l}$ be the total number of
stubs in the first group, and $k_1,\dots,k_{\sN-j}$ and $B_{\sN-j}$ the nodes and
number of stubs in the second group. We remark that $A_j+B_{\sN-j}=L_{\sN}$,
and since the groups are not connected, both $A_j$ and $B_{\sN-j}$ are even. The
$\prob_{\sN}$-probability that the groups are not connected is then,
for each fixed choice $i_1,\dots,i_j$, equal to
    \begin{equation}
    \label{dnz21}
    \produ\limits_{n=0}^{\frac{A_j}{2}-1}\frac{A_j-2n-1}{L_{\sN}-2n-1}
    =\frac{\produ\limits_{n=0}^{\frac{A_j}{2}-1}(A_j-2n-1)}
    {\produ\limits_{n=0}^{\frac{A_j}{2}-1}(L_{\sN}-2n-1)}
    =\produ\limits_{m=0}^{\frac{A_j}{2}-1}\frac{2m+1}{L_{\sN}-2m-1}.
    \end{equation}
By symmetry between $A_j$ and $B_{\sN-j}$,
we also obtain that this probability is equal to
    \eq
    \label{dnz21b}
    \produ\limits_{m=0}^{\frac{B_{N-j}}{2}-1}\frac{2m+1}{L_{\sN}-2m-1}.
    \en
Observe that for integers $j\geq 0$, the map
    \begin{equation}
    \label{dnz22}
    j\mapsto\produ\limits_{m=0}^j\frac{2m+1}{L_{\sN}-2m-1}\qquad
    \mbox{ is decreasing for }\qquad j\le\frac{L_{\sN}}{4}-\frac{1}{2}.
    \end{equation}
Suppose that $A_j\le L_{\sN}/2-1$. Then we use (\ref{dnz21}). Due to $\prob(D_1\ge r)=1$,
and since $A_j$ is even
we have $\lceil jr/2\rceil \le A_j/2\le L_{\sN}/4-1/2$ a.s., and hence,
by~(\ref{dnz22}), the final
expression in~(\ref{dnz21}) is at most
    $$
    \produ\limits_{m=0}^{\lceil jr/2\rceil -1}\frac{2m+1}{L_{\sN}-2m-1}.
    $$
Suppose that
$A_j\ge L_{\sN}/2$. Since $A_j+B_{\sN-j}=L_{\sN}$ and $\prob(D_1\ge r)=1$,
we have then that $\lceil (N-j)r/2\rceil-1 \le B_{\sN-j}/2-1\le L_{\sN}/4-1/2$
a.s., and we estimate~(\ref{dnz21b}),
by~(\ref{dnz22}), by
    $$
    \produ\limits_{m=0}^{\lceil (N-j)r/2\rceil -1}\frac{2m+1}{L_{\sN}-2m-1}.
    $$
Hence, the $\prob_{\sN}$-probability that the two groups of nodes are not connected
is at most
    \begin{equation}
    \label{dnz26}
    \begin{array}{rl}
    \produ\limits_{m=0}^{\lceil jr/2\rceil-1}\frac{2m+1}{L_{\sN}-2m-1}
    \indic[A_j\le L_{\sN}/2-1]&+
    \produ\limits_{m=0}^{\lceil (N-j)r/2\rceil-1}\frac{2m+1}{L_{\sN}-2m-1}
    \indic[A_j\geq L_{\sN}/2]\\
    \le\produ\limits_{m=0}^{\lceil jr/2\rceil-1}\frac{2m+1}{L_{\sN}-2m-1}&+
    \produ\limits_{m=0}^{\lceil (N-j)r/2\rceil-1}\frac{2m+1}{L_{\sN}-2m-1}.
    \end{array}
    \end{equation}
For $1\le j\le N-1$, we have at most $N \choose j$ ways to
choose $j$ nodes $i_1, \ldots, i_j$. Hence, by Boole's inequality,
    \begin{equation}
    \label{dnz20}
    \begin{array}{rl}
    \prob_{\sN}\left(\CN\ge s\right)\le
    &\sumu\limits_{j=s}^{N-s} \frac{N!}{j!(N-j)!}\left(
    \produ\limits_{m=0}^{\lceil jr/2\rceil-1}\frac{2m+1}{L_{\sN}-2m-1}+
    \produ\limits_{m=0}^{\lceil (N-j)r/2\rceil-1}\frac{2m+1}{L_{\sN}-2m-1}\right)\\
    &\qquad = 2\sumu\limits_{j=s}^{N-s} \left(\produ\limits_{m=\lceil
    jr/2\rceil}^{j-1}\frac{N-m}{m+1}\right)\left(\produ\limits_{m=0}^{\lceil jr/2\rceil-1}
    \frac{(N-m)(2m+1)}{(m+1)(L_{\sN}-2m-1)}\right),
    \end{array}
    \end{equation}
where by convention the product of the empty set equals $1$,
and where we used symmetry (between $s$ and $N-s$) together with the identity
$$
\frac{N!}{j!(N-j)!}=\prod_{m=0}^{j-1} \frac{N-m}{m+1}.
$$
For the remaining part of the proof we will make use of the following lemma:
\begin{lemma}
For any $1\le k\le N-1$,
\begin{equation}
\label{dnz24} \produ\limits_{m=0}^{k-1}
\frac{(N-m)(2m+1)}{(m+1)(2N-2m-1)}\le1.
\end{equation}
\end{lemma}

\proof Define, for $0\leq m\leq N-1$,
    $$
    h(m)=\frac{(N-m)(2m+1)}{(m+1)(2N-2m-1)}.
    $$
Then
    $$
    \begin{array}{rl}
    h(m)\le1,&\mbox{ if }m\le(N-1)/2,\\
    h(m)h(N-m-1)=1,&\mbox{ for all }0\le m\leq N-1,\\
    h\left((N-1)/2\right)=1,&\mbox{ if }N\mbox{ is odd}.
    \end{array}
    $$
Hence,~(\ref{dnz24}) is trivial for $k\le(N-1)/2+1$. If
$k>(N-1)/2+1$ then $N-k<(N-1)/2$, and
$$
\left(\produ\limits_{m=0}^{k-1} h(m)\right)^2
\le\left(\produ\limits_{m=N-k}^{k-1}h(m)\right)^2
=\produ\limits_{m=N-k}^{k-1}h(m)h(N-m-1)=1.$$
Thus we have~(\ref{dnz24}) for all $1\le k\le N-1$.
\qed
\vskip 0.5cm

\noindent
We now finish the proof of Proposition \ref{dnzLeTh1} when $r=2$.
In this case, due to~(\ref{dnz24}), the right side of (\ref{dnz20}), with $r=2$, is at most
    \eq
    2\sumu\limits_{j=s}^{N-s} \produ\limits_{m=0}^{j-1}
        \frac{(N-m)(2m+1)}{(m+1)(L_{\sN}-2m-1)}\leq
        2\sumu\limits_{j=s}^{N-s} \produ\limits_{m=0}^{j-1}
    \frac{2N-2m-1}{L_{\sN}-2m-1} \le 2\sumu\limits_{j=s}^{N-s}
    \left(\frac{2N}{L_{\sN}}\right)^j,
    \en
since $L_{\sN}\geq 2N,$ a.s. This completes the proof when $r= 2$.

When $r=1$, the right side of (\ref{dnz20}) equals
    \eq
    2\sumu\limits_{j=s}^{N-s}\left(\produ\limits_{m=\lceil
    j/2\rceil}^{j-1}\frac{N-m}{m+1}\right)\left(\produ\limits_{m=0}^{\lceil
    j/2\rceil-1} \frac{(N-m)(2m+1)}{(m+1)(L_{\sN}-2m-1)}\right),
    \en
which,
due to~(\ref{dnz24}), is at most
    \eq
    2\sumu\limits_{j=s}^{N-s}\left(\produ\limits_{m=\lceil
    j/2\rceil}^{j-1}\frac{N-m}{m+1}\right)\left(\produ\limits_{m=0}^{\lceil
    j/2\rceil-1} \frac{2N-2m-1}{L_{\sN}-2m-1}\right)\le
    2\sumu\limits_{j=s}^{N-s}N^{\lfloor j/2\rfloor}\left(\frac{2N}{L_{\sN}}\right)^{\lceil j/2\rceil}
    \le2\sumu\limits_{j=s}^{N-s}\left(\frac{2N^2}{L_{\sN}}\right)^{\lceil
    j/2\rceil}.
    \nn
    \en
This completes the proof of Proposition \ref{dnzLeTh1}.
\qed


\section{On the connected component sizes}
\label{sec-conncom}
In this section, we investigate the largest connected component in more detail and prove
Theorem \ref{thm-clusterP}.
We start with some definitions. For $\delta, \vep>0$, we define
$\gamma_{\sN}=\gamma_{\sN}(\delta, \vep)$ by
    \eq
    \label{gammaNdef}
    \gamma_{\sN}=\frac{1+\delta}{\log{\mu_{\sN}}-\log{2}-\frac{\vep}{1-2\vep}} \log{N},
    \en
where as before $\mu_{\sN}=L_{\sN}/N$. We also define a deterministic version of
$\gamma_{\sN}$ in the following way:
    \eq
    \label{gammabarNdef}
    \bar{\gamma}_{\sN} =\frac{1+\delta}{\log{\underline{\mu}_{\sN}}-\log{2}-\frac{\vep}{1-2\vep}}\log{N},
    \en
where $\underline{\mu}_{\sN}$ is a deterministic sequence for which
    \eq
    \label{mubardef}
    \prob(\mu_{\sN}\geq \underline{\mu}_{\sN})=1-o(1), \quad N\to \infty.
    \en
We start by formulating a version of Theorem \ref{thm-clusterP} that is valid
under the $\prob_{\sN}$-probability. Before stating this theorem, we need a
number of assumptions. Define
    \eq
    \label{ass3}
    q_{\sN} = \frac1N\sum_{i=1}^N \prob_{\sN}(|\Ccal_{i}|\geq \gamma_{\sN}),
    \en
where $\Ccal_i$ is the connected component that contains $i$ and
$|\Ccal|$ denotes the number of nodes in $\Ccal\subseteq\{1, \ldots, N\}$.
We assume that
    \eq
    \label{ass1}
    \frac{1}{N(N-1)}\sum_{i\neq j}\prob_{\sN}(i,j\text{ connected})=q_{\sN}^2(1+o(1)).
    \en
Note that (\ref{ass1}) is an assumption involving the $\prob_{\sN}$-probability.
We can interpret (\ref{ass1}) as saying that, under $\prob_{\sN}$,
a large proportion of nodes $i,j$ for which the connected component
consists of more than $\gamma_{\sN}$ nodes, are connected.

Before we proceed with the preliminaries of the proof of Theorem \ref{thm-clusterP} we give an
outline of this proof. Denote  for $\tau>2$\ by $q$  the survival probability of
the branching process $\{ {\cal Z}_l\}$ and set $q=1$ for $\tau\in (1,2)$. We will show,
using the coupling in \cite{HHV05} that
\begin{enumerate}
\item[(i)]
$q_{\sN}\to q$, see Lemma \ref{lem-qNasy}
 below,

\item[(ii)]
$\mbox{Var}_{\sN}(X_{\sN})=o(N^2)$, where $ X_{\sN}=\sum_{i=1}^N I[|\Ccal_i|\geq \bar{\gamma}_{\sN}]$,
see Lemma \ref{lem-varbdP} below.
\end{enumerate}
Having verified these two items, we can apply Proposition \ref{prop-red} below, because
$L_{\sN}\ge 2N$ follows from $\mu >2$ when $\tau>2$ holds, and is immediate for
$\tau\in (1,2)$. In the proof of Proposition \ref{prop-red} we will order the connected components
according to their size: $|\Ccal_{\sss (1)}|\geq |\Ccal_{\sss (2)}|\geq \ldots$, and will prove that
$$
\prob(|\Ccal_{\sss (2)}|\ge \bar{\gamma}_{\sN})=o(1),
$$
in two steps. In the proof of these two steps the statement of Lemma 3.4 below:
{\it the probability that
there exists a connected component with at most $\eta N$
nodes, and in between $\vep L_{\sN}$ and $(1-\vep)L_{\sN}$ stubs, is
exponentially small in $N$} plays a prominent role. We then finish with a proof that {\bf whp},
$$
 X_{\sN}=\sum_{i=1}^N I[|\Ccal_i|\geq \bar{\gamma}_{\sN}]
=\sum_{l}|\Ccal_{\sss(l)}|I[|\Ccal_{\sss(l)}|\geq \bar{\gamma}_{\sN}]
    =|\Ccal_{\sss (1)}|I[|\Ccal_{\sss (1)}|\geq \bar{\gamma}_{\sN}]
        =|\Ccal_{\sss (1)}|,
$$
and reach the conclusion that $X_{\sN}=|\Ccal_{\sss (1)}|$ is of order $Nq_{\sN}$ by showing that
$\prob(|X_{\sN}-Nq_{\sN}|> \omega_{\sN}\mbox{Var}(X_{\sN}))=o(1)$, for each sequence $\omega_{\sN}\to \infty$.

We now turn to the preliminaries of the proof of Theorem \ref{thm-clusterP}.
\begin{theorem}
    \label{thm-cluster}
    Assume that: (i) $L_{\sN}\geq 2N$, (ii) Relation (\ref{ass1}) holds, and  (iii) $q_{\sN}\geq \eps$ as $N\rightarrow \infty$,
    for some $\eps>0$. Then, {\bf whp}, under $\prob_{\sN}$, the
    largest connected component in $G$ has
    $q_{\sN}N(1+o(1))$ nodes, and all other connected components have at most
    $\gamma_{\sN}$ nodes. Moreover, {\bf whp},  under $\prob_{\sN}$, the largest connected component has
    in between $q_{\sN}N(1 \pm \omega_{\sN}\sqrt{\frac{\gamma_{\sN}}{L_{\sN}}})$
    nodes for any $\omega_{\sN}\rightarrow \infty$.
\end{theorem}

\begin{remark}
Observe that besides the results on the sizes of the components, Theorem \ref{thm-cluster}
also includes a bound on the fluctuation of the size of the largest component.
\end{remark}

The remainder of this section is organized as follows.
In Section \ref{sec-proofThmcs},
we prove Theorem \ref{thm-cluster}. In Section \ref{sec-pfThmclusterP}, we
use a modification of the proof of Theorem \ref{thm-cluster}
to prove Theorem \ref{thm-clusterP}.

\subsection{Connected components under $\prob_{\sN}$}
\label{sec-proofThmcs}
We start with a proposition that shows that the connected components,
measured in terms of their number of {\it edges}, are either quite small, i.e.,
less than $\gamma_{\sN},$ or very large, i.e., a positive fraction of the
total number of edges.

\begin{prop}
\label{thm-cs} Fix $\delta>0$, assume that
$\mu_{\sN}>2$, and let $0<\vep<\frac 1{10}$ be such that
$\log{\mu_{\sN}}-\log{2}-\frac{\vep}{1-2\vep}>0$.
Then, the $\prob_{\sN}$-probability that there
exists a connected component with in between
$\gamma_{\sN}=\gamma_{\sN}(\vep,\delta)$ and $\vep L_{\sN}$ edges
is bounded by $O(N^{-\delta})$.\\
Consequently, when $\mu>2$ or $\mu=\infty$,
the $\prob$-probability that there
exists a connected component with in between
$\bar\gamma_{\sN}$ and $\vep L_{\sN}$ edges
converges to 0, as $N\to\infty$.
\end{prop}

Of course, we expect that there is a unique
such large connected component, and this is what we will prove later on.
Note that when $\tau\in (1,2)$, then, with large probability, $\mu_{\sN}=L_{\sN}/N\geq N^{\eta},$
for some $\eta>0$. In this case, we even have that $\gamma_{\sN}$ is uniformly bounded in
$N$, and thus, the connected components that do not contain a positive fraction
of the edges are uniformly bounded in their number of edges.
\vskip0.5cm

\noindent
{\bf Proof of Proposition \ref{thm-cs}.}
We adapt the proof of Proposition \ref{dnzLeTh1}. Denote by $k$ the
number of edges in the connected component with in between
$\gamma_{\sN}$ and $\vep L_{\sN}$ edges. Then, we must have that all the
$2k$ stubs are connected to each other, i.e., they are not connected to
stubs not in the $k$ edges. This probability is bounded by
    $$
    \produ\limits_{n=0}^{k-1}\frac{2k-2n-1}{L_{\sN}-2n-1}
    =\produ\limits_{m=0}^{k-1}\frac{2m+1}{L_{\sN}-2m-1},
    $$
ignoring the fact that the component needs to be connected.
We first prove the statement for $k\leq (\frac{N}{2}-1)\wedge \vep
L_{\sN}$, and in a second step prove the statement for
$(\frac{N}{2}-1)\wedge \vep L_{\sN} < k \leq \vep L_{\sN}$. We
abbreviate $R_{\sN}=(\frac{N}{2}-1)\wedge \vep L_{\sN}$.

Denote the number of nodes in the connected component by $l$.
Note that when a connected component consists of $k$ edges, then
$l\leq k+1$. Therefore,
the total number of ways in which we can choose these $l$ nodes is
at most ${{N}\choose{l}}\leq {{N}\choose{k+1}},$
when $k\leq \frac{N}{2}-1$. Thus, the $\prob_{\sN}-$probability
that there exists a component with in between
$\gamma_{\sN}$ and $R_{\sN}$ edges is bounded by
    \eqalign
   \sum_{k=\gamma_{\sN}}^{R_{\sN}}
   {{N}\choose{k+1}}\produ\limits_{m=0}^{k-1}\frac{2m+1}{L_{\sN}-2m-1}
   &\leq  \sum_{k=\gamma_{\sN}}^{R_{\sN}}\frac{N!}{(N-k-1)!} 2^k
   \produ\limits_{m=0}^{k-1} (L_{\sN}-2m-1)^{-1}\nn\\
   &\leq N\sum_{k=\gamma_{\sN}}^{R_{\sN}} \big(\frac{2}{\mu_{\sN}}\big)^k
   \produ\limits_{m=0}^{k-1} (1+\frac{2m+1}{L_{\sN}-2m-1}).
   \enalign
Next, we use that $1+x\leq e^x$ for $x\geq 0$,
to obtain that the $\prob_{\sN}-$probability that there exists a component with in between
$\gamma_{\sN}$ and $R_{\sN}$ edges is bounded by
    \eqalign
   N\sum_{k=\gamma_{\sN}}^{R_{\sN}} \big(\frac{2}{\mu_{\sN}}\big)^k
   \exp\left\{\sum_{m=0}^{k-1} \frac{2m+1}{L_{\sN}-2m-1}\right\}
   &\leq N\sum_{k=\gamma_{\sN}}^{R_{\sN}} \big(\frac{2}{\mu_{\sN}}\big)^k
   e^{\frac{k^2}{L_{\sN}-2\vep L_{\sN}}}\nn\\
   &\leq N\sum_{k=\gamma_{\sN}}^{R_{\sN}} \big(\frac{2}{\mu_{\sN}}\big)^k e^{\frac{k\vep}{1-2\vep}}
   \leq \eta^{-1} N\Big(\frac{2}{\mu_{\sN}}e^{\frac{\vep}{1-2\vep}}\Big)^{\gamma_{\sN}},
   \label{finalbd}
   \enalign
provided that $\frac{2}{\mu_{\sN}}e^{\frac{\vep}{1-2\vep}}\le 1-\eta$.
The right side of (\ref{finalbd}) is bounded by $\eta^{-1} N^{-\delta}$
for the choice of $\gamma_{\sN}$ in (\ref{gammaNdef}).

We complete the proof by dealing with the case that $R_{\sN}<k\leq\vep L_{\sN}$.
In this case, we must have that $R_{\sN}=\frac{N}{2}-1$, otherwise there is nothing to prove, so that $k\geq \frac{N}{2}$.
Then, we bound the total number of ways
in which we can choose the $l\leq k+1$ nodes by $2^N$. Since $2^{N}
\leq 2^{2k},$ for all $k\geq N/2$, we arrive at the fact that the probability
that there exists a connected component with in between $\frac{N}{2}$ and $\vep L_{\sN}$
edges is bounded by
    \eqalign
   \sum_{k=\frac{N}{2}}^{\vep L_{\sN}}
    2^N \produ\limits_{m=0}^{k-1}\frac{2m+1}{L_{\sN}-2m-1}
   &\leq  \sum_{k=\frac{N}{2}}^{\vep L_{\sN}} 2^{2k}\big(\frac{2\vep}{1-2\vep}\big)^k\leq
   \frac{1-2\vep}{1-10\vep}\left(
\frac{8\vep}{1-2\vep}
   \right)^{N/2},
   \enalign
which is exponentially small in $N$ when $\vep<\frac 1{10}$. Thus, this probability is certainly bounded
above by $N^{-\delta}$.

For the bound on the $\prob$-probability that there
exists a connected component with in between
$\bar\gamma_{\sN}$ and $\vep L_{\sN}$ edges, we denote by
$F(k,l)$ the event that there exists a connected component with
in between $k$ and $l$ edges. Then we can bound
    \eq
    \prob(F(\bar\gamma_{\sN},\vep L_{\sN}))\leq
    \prob(\mu_{\sN}< \underline{\mu}_{\sN})+
    \expec\big[I[\mu_{\sN}\geq \underline{\mu}_{\sN}]\prob_{\sN}(F(\gamma_{\sN},\vep L_{\sN}))\big],
    \en
where we use that $\bar\gamma_{\sN}\geq \gamma_{\sN}$ when
$\mu_{\sN}\geq \underline{\mu}_{\sN}$, choosing $ \underline{\mu}_{\sN}=(\mu+2)/2$ for $\mu<\infty$ and
$ \underline{\mu}_{\sN}=3$ for $\mu=\infty$.
The first term is $o(1)$ by (\ref{mubardef}), the second term is small by
the estimate $\prob_{\sN}(F(\gamma_{\sN},\vep L_{\sN}))\leq N^{-\delta}$
proved above.
\qed
\vskip0.5cm

We next present a lemma that will be used in the proof of
Theorem \ref{thm-cluster}.

\begin{lemma}
\label{lem-smnod, largestubs}
Fix $\eps>0$ and $0<\eta<\vep$ sufficiently small.
Then, when $L_{\sN}\geq 2N$, the $\prob_{\sN}$-probability that
there exists a connected component with at most $\eta N$
nodes, and in between $\vep L_{\sN}$ and $(1-\vep)L_{\sN}$ stubs, is
exponentially small in $N$. Consequently, the same estimate is true for
the $\prob$-probability of this event provided that $\prob(L_{\sN}< 2N)$
is exponentially small.
\end{lemma}

\proof
Take  $0< \eta<\vep$. Again denote by $k$ the
number of edges in an arbitrary connected component satisfying
$\frac{\vep}{2} L_{\sN} \leq k \leq
\frac{(1-\vep)}{2}L_{\sN}$. Then, we must have that all the
$2k$ stubs are connected to each other, i.e., they are not connected to
stubs not in the $k$ edges. The $\prob_{\sN}$-probability of this event
is bounded by
    $$
    \produ\limits_{n=0}^{k-1}\frac{2k-2n-1}{L_{\sN}-2n-1}
    \leq \produ\limits_{n=0}^{k-1}\frac{2k-2n}{L_{\sN}-2n}=\produ\limits_{n=0}^{k-1}\frac{k-n}{\frac{L_{\sN}}{2}-n}
    ={{\frac{L_{\sN}}{2}}\choose{k}}^{-1}.
    $$
We next use that the number ways of choosing
at most $\eta N$  nodes, with $\eta<\frac12$, is bounded from above by
    \eq
    \sum_{j=0}^{\eta N} {{N}\choose{j}}
    \leq (\eta N+1){{N}\choose{\eta N}}<
    N{{N}\choose{\eta N}}.
    \en
Therefore, the $\prob_{\sN}$-probability that
there exists a connected component with in between $\vep L_{\sN}$ and $(1-\vep)L_{\sN}$ stubs
and at most $\eta N$ nodes is bounded by
\begin{eqnarray}
&&N{{N}\choose{\eta N}} \sum_{k=\frac{\vep}{2} L_{\sN}}^{\frac{(1-\vep)}{2}L_{\sN}}
    {{\frac{L_{\sN}}{2}}\choose{k}}^{-1}\nn
    \leq NL_{\sN}{{N}\choose{\eta N}} {{\frac{L_{\sN}}{2}}\choose{\frac{\vep}{2}L_{\sN}}}^{-1}\\
&&\quad
= NL_{\sN} \exp\{c_\eta N(1+o(1))\}\exp\{-c_{\vep}\frac{L_{\sN}}{2}(1+o(1))\},
\end{eqnarray}
where we have bound the sum by the number terms times
the largest summand, and where we have used that for $\eta$  small, ${{N}\choose{\eta N}}=e^{c_{\eta}N(1+o(1))}$,
where $c_{\eta}\downarrow 0$ as $\eta\downarrow 0$. Therefore, using that $L_{\sN}\geq 2N$,
it suffices to take $\eta>0$ so small that
$c_\eta N < (c_{\vep}-\delta)\frac{L_{\sN}}{2}$, for some $\delta>0$ sufficiently small,
to see that this probability is exponentially small in $N$.

For the statement involving the unconditional probability we denote by $G_{\vep,\eta}$ the event that
there exists a connected component with at most $\eta N$
nodes, and in between $\vep L_{\sN}$ and $(1-\vep)L_{\sN}$ stubs. Then
\eq
\label{toevoeging1}
\prob(G_{\vep,\eta})=\expec[\prob_{\sN}(G_{\vep,\eta})]
\leq
\prob(L_{\sN}<2N)+\expec[\prob_{\sN}(G_{\vep,\eta}, L_{\sN}\geq 2N)]
\en
and both terms are exponentially small.
This completes the proof of Lemma \ref{lem-smnod, largestubs}.
\qed
\vskip0.5cm

\noindent
We are now ready for the proof of Theorem \ref{thm-cluster}:\\
{\bf Proof of Theorem \ref{thm-cluster}.} Take  $\vep,\delta>0$ and fix $\gamma_{\sN}$
as in (\ref{gammaNdef}). Recall that $\Ccal_i$ denotes the connected
component that $i$ belongs to.
We define the random variable $X_{\sN}$ by
    \eq
    X_{\sN}=\sum_{i=1}^N I[|\Ccal_i|\geq \gamma_{\sN}],
    \en
so that $X_{\sN}$ equals the total number of nodes in connected components of size at least
$\gamma_{\sN}$. By (\ref{ass3}),
    \eq
    \label{expXN}
    \expec_{\sN}[X_{\sN}] = \sum_{i=1}^N \prob_{\sN}(|\Ccal_i|\geq \gamma_{\sN})
    =N q_{\sN},
    \en
where $\expec_{\sN}$ is the expected value under\ $\prob_{\sN}$.
We first prove that the variance of $X_{\sN}$ under the law $\prob_{\sN}$
is small, so that $X_{\sN}$ is with high probability
close to $Nq_{\sN}$:

\begin{lemma}
\label{lem-varbd}
With probability 1,
\eq
{\rm Var}_{\sN}(X_{\sN})=Nq_{\sN}(1-q_{\sN})+O(\frac{\gamma_{\sN}^2N^2}{L_{\sN}}),
\en
where ${\rm Var}_{\sN}$ denotes the variance under $\prob_{\sN}$.
\end{lemma}

\proof
Without explicit mentioning all statements in the proof hold with probability 1.
We first note that
$\text{Var}_{\sN}(X_{\sN})=\text{Var}_{\sN}(N-X_{\sN})
=\text{Var}_{\sN}(Y_{\sN})$, where
    \eq
    Y_{\sN}=\sum_{i=1}^N I[|\Ccal_i|< \gamma_{\sN}].
    \en
Therefore,
    \eqalign
    \text{Var}_{\sN}(Y_{\sN})&=\sum_{i,j}
    \prob_{\sN}(|\Ccal_i|<\gamma_{\sN},|\Ccal_j| <
    \gamma_{\sN})-[N(1-q_{\sN})]^2
    \label{YN1}\\
    &= \sum_{i\neq j}
    \prob_{\sN}(|\Ccal_i|<\gamma_{\sN},|\Ccal_j| <
    \gamma_{\sN})+N(1-q_{\sN})-N^2(1-q_{\sN})^2.\nn
    \enalign
For the first term we use the coupling
in \cite[Proof of Lemma A.2.2]{HHV05}, with $N^{\frac12-\eta}$  replaced by $\gamma_{\sN}$, to obtain that
    \eq
    \label{couplingZhat}
    \prob_{\sN}(|\Ccal_i|<\gamma_{\sN})=
    \prob_{\sN}(\sum_{l} \hat Z_{l}^{\smallsup{i,N}}<\gamma_{\sN})+O(\frac{\gamma_{\sN}^2}{L_{\sN}}),
    \en
where $\{\hat Z_{l}^{\smallsup{i,N}}\}_{l\geq 1}$ is a branching process with offspring
distribution
     \eq
     g_n^{\smallsup{N}}=\frac{n+1}{L_{\sN}}\sum_{j=1}^N I[D_j=n+1],\qquad n\ge 0,
     \en
and with $\hat Z_{1}^{\smallsup{i,N}}=D_i$,
the degree of node $i$. The coupling is described in full detail
in \cite[Section 3]{HHV05}, whereas the bound in (\ref{couplingZhat})
follows from the proof of \cite[Lemma A.2.2]{HHV05}, which holds under the
$\prob_{\sN}$-probability and is therefore true for any degree sequence, and hence
in particular for each $\tau>1$.

Moreover, for $i\neq j$, it is described in \cite[Section 3]{HHV05} that
we can couple $|\Ccal_i|$ and $|\Ccal_j|$ simultaneously to two
{\it independent} branching processes to obtain
    \eq
    \prob_{\sN}(|\Ccal_i|<\gamma_{\sN},|\Ccal_j| <
    \gamma_{\sN})=\prob_{\sN}(\sum_{l} \hat Z_{l}^{\smallsup{i,N}}<\gamma_{\sN})
    \prob_{\sN}(\sum_{l} \hat Z_{l}^{\smallsup{j,N}}<\gamma_{\sN})+O(\frac{\gamma_{\sN}^2}{L_{\sN}}).
    \en
Therefore,
    \eqalign
    \sum_{i\neq j} \prob_{\sN}(|\Ccal_i|<\gamma_{\sN},|\Ccal_j| <
    \gamma_{\sN}) &= \sum_{i\neq j}\big[\prob_{\sN}(\sum_{l} \hat Z_{l}^{\smallsup{i,N}}<\gamma_{\sN})
    \prob_{\sN}(\sum_{l} \hat Z_{l}^{\smallsup{j,N}}<\gamma_{\sN})+O(\frac{\gamma_{\sN}^2}{L_{\sN}})\big]\nn\\
    &= \Big(\sum_{i=1}^N\prob_{\sN}(\sum_{l} \hat Z_{l}^{\smallsup{i,N}}<\gamma_{\sN})\Big)^2
    - (1-q_{\sN})^2 N+O(\frac{\gamma_{\sN}^2N^2}{L_{\sN}})\nn\\
    &=(N^2-N)(1-q_{\sN})^2  +O(\frac{\gamma_{\sN}^2N^2}{L_{\sN}}),
    \label{YN2}
    \enalign
using (\ref{ass3}) and (\ref{couplingZhat}). So,
substituting (\ref{YN2}) into (\ref{YN1}),
    \begin{equation}
    \text{Var}_{\sN}(X_{\sN})=\text{Var}_{\sN}(Y_{\sN})=
     Nq_{\sN}(1-q_{\sN})+ O(\frac{\gamma_{\sN}^2N^2}{L_{\sN}}).
    \end{equation}
\qed

\noindent  We continue with the proof of Theorem
\ref{thm-cluster}, which is a consequence of the following
proposition. This proposition will also be used to prove Theorem
\ref{thm-clusterP} below. In its statement, we let $\mathbb{Q}$ be
a probability distribution, which we will take to be $\prob_{\sN}$
in the proof of Theorem \ref{thm-cluster} and $\prob$ in the proof
of Theorem \ref{thm-clusterP}. Let $\gamma_{\sN}^*=\gamma_{\sN}$
when $\mathbb{Q}=\prob_{\sN}$ and $\gamma_{\sN}^*=\bar
\gamma_{\sN}$ when $\mathbb{Q}=\prob$,
see (\ref{gammaNdef})
 and (\ref{gammabarNdef}) for the definitions of $\gamma_{\sN}$ an $\bar
\gamma_{\sN}$.
Furthermore, we take $X_{\sN}=\sum_{i=1}^N I[|\Ccal_i|\geq
\gamma^*_{\sN}]$ and define
\eq
    \label{ass4}
    q^*_{\sN} = \frac1N\sum_{i=1}^N \Q(|\Ccal_{i}|\geq \gamma^*_{\sN}).
    \en

\begin{prop}
\label{prop-red}Let $\Q=\prob$ or $\Q=\prob_{\sN}$.
Suppose that (i) $L_N\geq 2N$, (ii) $\text{Var}_{\sss \Q}(X_{\sN})\leq B_{\sN}=o(N^2)$, and
    \eq
    (iii)\quad
    \expec_{\sss \Q}[X_{\sN}]=Nq^*_{\sN},\qquad     \sum_{i,j} \mathbb{Q}(i,j \text{ connected})
    =(Nq^*_{\sN})^2(1+o(1)),
    \label{ass3rep}
    \en
where $q^*_{\sN}\geq \eps$ for some $\eps>0$, as $N\to \infty$. Then,
\begin{itemize}
\item[{\rm (i)}]
\whp the second largest component has at most $\gamma^*_{\sN}$ nodes;
\item[{\rm (ii)}]
\whp the largest connected component has
in between $Nq^*_{\sN} \pm \omega_{\sN}\sqrt{B_{\sN}}$ nodes
for any $\omega_{\sN}\rightarrow \infty$,
such that $\omega_{\sN}\sqrt{B_{\sN}}=o(N)$.
\end{itemize}
\end{prop}
To prove Theorem \ref{thm-cluster}, we use the above with $\Q=\prob_{\sN}$
and $B_{\sN}=C\frac{\gamma_{\sN}^2N^2}{L_{\sN}}$.

\proof
We define the event
    \eq
    E_{\sN} = \{|X_{\sN}-Nq^*_{\sN}|\leq \omega_{\sN}\sqrt{B_{\sN}}\}.
    \en
Then, by the Chebycheff inequality,
    \eq
    \label{ENcomplbd}
    \Q(E_{\sN}^c)\leq \Big(\omega_{\sN}\sqrt{B_{\sN}}\Big)^{-2}
    \text{Var}_{\sss \Q}(X_{\sN})\leq \omega_{\sN}^{-2}=o(1).
    \en

We write $\Ccal_{\sss(1)}, \Ccal_{\sss(2)}, \ldots$
for the connected components ordered according to their sizes, so
that $|\Ccal_{\sss(1)}|\geq |\Ccal_{\sss(2)}|\geq \ldots$ and
$\Ccal_{\sss(i)}$ and  $\Ccal_{\sss(j)}$ are disjoint for $i\neq j$.
Then we clearly have that
    \eq
    \label{rewC(l)s}
    \sum_{i,j} \Q(i,j\text{ connected})=\sum_{i,j} \Q\Big(\bigcup_{l} \{i,j\in \Ccal_{\sss(l)}\}\Big)=
    \sum_{l}\sum_{i,j} \Q(i,j\in \Ccal_{\sss(l)})
    =\sum_{l} \expec_{\sss \Q}[|\Ccal_{\sss(l)}|^2].
    \en
Combining with (\ref{ass3rep}) we get,
    \eq
    \label{asscons}
    \sum_{l} \expec_{\sss \Q}[|\Ccal_{\sss(l)}|^2]=(Nq^*_{\sN})^2 (1+o(1)).
    \en
Furthermore,
    \eqalign
    \sum_{l} \expec_{\sss \Q}\big[|\Ccal_{\sss(l)}|^2I[|\Ccal_{\sss(l)}|< \gamma^*_{\sN}]\big]
    &\leq \gamma_{\sN}^* \sum_{l} \expec_{\sss \Q}\big[|\Ccal_{\sss(l)}|I[|\Ccal_{\sss(l)}|<
    \gamma^*_{\sN}]\big]\leq \gamma^*_{\sN} N.
    \enalign
Therefore, since $\gamma^*_{\sN}=O(\log{N})=o(N)$ and $q^*_{\sN}\geq \eps$,
we obtain that
    \eq
    \sum_{l} \expec_{\sss \Q}\big[|\Ccal_{\sss(l)}|^2I[|\Ccal_{\sss(l)}|\geq \gamma^*_{\sN}]\big]
    =(Nq^*_{\sN})^2 (1+o(1)).
    \en
By (\ref{ENcomplbd}), we thus also have that
    \eq
    \label{suml2}
    \expec_{\sss \Q}\big[\sum_{l}|\Ccal_{\sss(l)}|^2
    I[|\Ccal_{\sss(l)}|\geq \gamma^*_{\sN}]I[E_{\sN}]\big]
    =(Nq^*_{\sN})^2 (1+o(1)).
    \en
We will now prove that
    \eq
    \label{seccompbd}
    \Q(|\Ccal_{\sss (2)}|\geq \gamma^*_{\sN})=o(1).
    \en
This proceeds in two key steps. We first show that
for some $\eta>0$ sufficiently small
    \eq
    \label{reduc1}
    \Q(|\Ccal_{\sss (2)}|\geq \gamma_{\sN}^*)=\Q(|\Ccal_{\sss (2)}|>\eta N)
    +o(1),
    \en
and then that the assumption that
    \eq
    \limsup_{N\rightarrow \infty}\Q(|\Ccal_{\sss (2)}|>\eta N)=\theta>0,
    \label{contrass}
    \en
leads to a contradiction. Together, this proves (\ref{seccompbd}).
We start by proving (\ref{reduc1}). We note that
we only need to prove that $\Q(|\Ccal_{\sss (2)}|\geq \gamma_{\sN}^*)$
is less than or equal to the right side of (\ref{reduc1}),
since the other bound is trivial (even with $o(1)$ replaced by 0).

To prove (\ref{reduc1}), we split for $i=1,2$,
    \eq
    \Q(|\Ccal_{\sss (i)}|\geq \gamma^*_{\sN})
    = \Q(|\Ccal_{\sss (i)}|\geq \gamma^*_{\sN}, |\Ccal_{\sss(i)}|_b\leq \vep L_{\sN})
    +\Q(|\Ccal_{\sss (i)}|\geq \gamma^*_{\sN}, |\Ccal_{\sss(i)}|_b>\vep L_{\sN}),
    \en
where $|\Ccal|_b$ denotes the number of edges in $\Ccal$. Since
$|\Ccal|_b\geq |\Ccal|-1$, for any connected component $\Ccal$, by Proposition \ref{thm-cs},
for any $\delta>0$, and for $i=1,2$, since $\gamma^*_{\sN}=\gamma_{\sN}$ or $\gamma^*_{\sN}
=\bar\gamma_{\sN}$, where \whp $\bar\gamma_{\sN}\geq \gamma_{\sN}$, we obtain
    \eq
    \Q(|\Ccal_{\sss (i)}|\geq\gamma^*_{\sN}, |\Ccal_{\sss(i)}|_b\leq \vep L_{\sN})
    \leq \Q(\gamma^*_{\sN}\leq |\Ccal_{\sss(i)}|_b\leq \vep L_{\sN})
    =o(1),
    \en
so that
    \eq
    \label{ger1606}
    \Q(|\Ccal_{\sss (2)}|\geq\gamma^*_{\sN})=\Q\big(|\Ccal_{\sss (2)}|\geq\gamma^*_{\sN},
    |\Ccal_{\sss(1)}|_b>\vep L_{\sN}, |\Ccal_{\sss(2)}|_b>\vep L_{\sN})+o(1).
    \en
By Lemma \ref{lem-smnod, largestubs},
and because $\{|\Ccal_{\sss(1)}|_b>\vep L_{\sN}\} \Rightarrow \{ |\Ccal_{\sss(2)}|_b<(1-\vep) L_{\sN}\}$,
we further have that for $\eta>0$ sufficiently small
    \eq
    \Q(|\Ccal_{\sss (2)}|\leq\eta N, |\Ccal_{\sss(1)}|_b>\vep L_{\sN},
    |\Ccal_{\sss(2)}|_b>\vep L_{\sN})=o(1),
    \en
Therefore, using (\ref{ger1606})
    \eq
    \Q(|\Ccal_{\sss (2)}|\geq \gamma^*_{\sN})\le
    \Q(|\Ccal_{\sss (2)}|>\eta N, |\Ccal_{\sss(1)}|_b>\vep L_{\sN},
    |\Ccal_{\sss(2)}|_b>\vep L_{\sN})+o(1)
    \leq \Q(|\Ccal_{\sss (2)}|>\eta N)
    +o(1).
    \en
This proves (\ref{reduc1}).

%
%

We next prove that (\ref{contrass}) is in contradiction
with (\ref{suml2}). Observe that
    \eqalign
    \label{onen}
X_{\sN}=\sum_{i=1}^N I[|\Ccal_i|\geq \gamma^*_{\sN}]
    =\sum_{l}|\Ccal_{\sss(l)}|I[|\Ccal_{\sss(l)}|\geq \gamma^*_{\sN}],
    \enalign
so that, on the event $E_{\sN}$, we have that
$\sum_{l}|\Ccal_{\sss(l)}|I[|\Ccal_{\sss(l)}|\geq \gamma^*_{\sN}]
=Nq^*_{\sN}(1+o(1))$. Using \eqref{onen} we can bound
    \eq
    \label{suml2b1}
    \sum_{l} |\Ccal_{\sss(l)}|^2
    I[|\Ccal_{\sss(l)}|\geq \gamma^*_{\sN}]
    \leq |\Ccal_{\sss(2)}|^2 +\big(\sum_{l\neq 2}|\Ccal_{\sss(l)}|I[|\Ccal_{\sss(l)}|\geq \gamma^*_{\sN}]\big)^2
    =|\Ccal_{\sss(2)}|^2 + (X_{\sN}-|\Ccal_{\sss(2)}|)^2.
    \en
We split the  expectation in (\ref{suml2}) by intersecting with
the event $\{|\Ccal_{\sss(2)}|>\eta N\}$ and its complement:
\begin{eqnarray}
\expec_{\sss \Q}\Big[\sum_{l}|\Ccal_{\sss(l)}|^2I[|\Ccal_{\sss(l)}|\geq \gamma^*_{\sN}]I[E_{\sN}]\Big]
&=&\expec_{\sss \Q}\Big[\sum_{l}|\Ccal_{\sss(l)}|^2
    I[|\Ccal_{\sss(l)}|\geq \gamma^*_{\sN}]I[E_{\sN}\cap \{|\Ccal_{\sss(2)}|>\eta N\}]\Big]
\label{hooghiem1}\\
&&\quad+\expec_{\sss \Q}\Big[\sum_{l}|\Ccal_{\sss(l)}|^2
    I[|\Ccal_{\sss(l)}|\geq \gamma^*_{\sN}]I[E_{\sN}\cap \{|\Ccal_{\sss(2)}|\leq\eta N\}]\Big].\nn
\end{eqnarray}
We next use a simple calculus argument. For $\eta N\leq x \leq y/2$, the function
$x\mapsto x^2+(y-x)^2$ is maximal for $x=\eta N$.
We apply the arising inequality to the right side (\ref{suml2b1}), with $x=|\Ccal_{\sss(2)}|$
and $y=X_{\sN}\ge |\Ccal_{\sss(1)}|+|\Ccal_{\sss(2)}|\geq 2|\Ccal_{\sss(2)}|=2x$, so that,
\begin{eqnarray}
    \label{suml2b}
&&\expec_{\sss \Q}\big[\sum_{l}|\Ccal_{\sss(l)}|^2
    I[|\Ccal_{\sss(l)}|\geq \gamma^*_{\sN}]I[E_{\sN}\cap \{|\Ccal_{\sss(2)}|>\eta N\}]\big]\nonumber\\
&&\qquad\qquad\leq\expec_{\sss \Q}\big[(|\Ccal_{\sss(2)}|^2+(X_{\sN}-|\Ccal_{\sss(2)}|)^2)
    I[E_{\sN}\cap \{|\Ccal_{\sss(2)}|>\eta N\}]\big]\nonumber\\
    &&\qquad\qquad\leq\expec_{\sss \Q}\big[(\eta^2N^2+(X_{\sN}-\eta^2N^2)
    I[E_{\sN}\cap \{|\Ccal_{\sss(2)}|>\eta N\}]\big]\nonumber\\
&&\qquad\qquad \leq \big(\eta^2+(q^*_{\sN}-\eta)^2\big)N^2\Q(|\Ccal_{\sss (2)}|>\eta N)(1+o(1)).
\end{eqnarray}
where we used in the last step that on $E_{\sN}$ we have $X_{\sN}=q^*_{\sN}N(1+o(1))$, because
$\omega_{\sN}\sqrt{B_{\sN}}=o(N)$.
On the other hand, we have, on the event $E_{\sN}$,
using again that $\omega_{\sN}\sqrt{B_{\sN}}=o(N)$,
    \eq
    \label{hooghiem2}
    \sum_{l} |\Ccal_{\sss(l)}|^2I[|\Ccal_{\sss(l)}|\geq \gamma^*_{\sN}]
    \leq \big(\sum_{l} |\Ccal_{\sss(l)}|I[|\Ccal_{\sss(l)}|
    \geq\gamma^*_{\sN}]\big)^2 = X_{\sN}^2= (Nq^*_{\sN})^2(1+o(1)),
    \en
implying that
\begin{equation}
\label{hooghiem3}
\expec_{\sss \Q}\big[\sum_{l}|\Ccal_{\sss(l)}|^2
    I[|\Ccal_{\sss(l)}|\geq \gamma^*_{\sN}]I[E_{\sN}\cap \{|\Ccal_{\sss(2)}|\leq\eta N\}]\big]
    \le \Q(\gamma^*_{\sN}\leq |\Ccal_{\sss (2)}|\leq \eta N)(Nq^*_{\sN})^2(1+o(1)).
\end{equation}
Together, (\ref{hooghiem1}), (\ref{suml2b}) and (\ref{hooghiem3}) yield
    \eqalign
    \label{hooghiem4}
    &\expec_{\sss \Q}\big[\sum_{l}|\Ccal_{\sss(l)}|^2
    I[|\Ccal_{\sss(l)}|\geq \gamma^*_{\sN}]I[E_{\sN}]\big]\\
    &\qquad \leq
    \Big[\big(\eta^2+(q^*_{\sN}-\eta)^2\big)N^2\Q(|\Ccal_{\sss (2)}|>\eta N)
    +(q^*_{\sN})^2N^2\Q(\gamma^*_{\sN}\leq |\Ccal_{\sss (2)}|\leq \eta N)\Big](1+o(1)),\nn
    \enalign
so that the assumption that
$\limsup_{N\rightarrow \infty} \Q   (|\Ccal_{\sss (2)}|>\eta N)= \theta>0$
is in contradiction with (\ref{suml2}), because assuming both (\ref{suml2}) and
(\ref{hooghiem4}) would imply that $\eta\geq \limsup  q_{\sN}^*=\vep$,
since for $0<\eta < q_{\sN}^*$ we have $\eta^2+( q_{\sN}^*-\eta)^2< (q_{\sN}^*)^2$. This proves that the assumption in
(\ref{contrass}) is false, and
we conclude that (\ref{seccompbd}) holds, which proves the claim for the second largest component.

We now prove that {\bf whp} the largest
component has size in between $Nq^*_{\sN} \pm \omega_{\sN}\sqrt{B_{\sN}}$
for any $\omega_{\sN}\rightarrow \infty$. For this, we note that on the event
that the second largest component
has size less than or equal to $\gamma^*_{\sN}$, we have (compare (\ref{onen})),
    \eq
    X_{\sN}=\sum_{l}|\Ccal_{\sss(l)}|I[|\Ccal_{\sss(l)}|\geq \gamma^*_{\sN}]
    =|\Ccal_{\sss (1)}|I[|\Ccal_{\sss (1)}|\geq \gamma_{\sN}^*]
        =|\Ccal_{\sss (1)}|.
    \en
By (\ref{ENcomplbd}) and (\ref{seccompbd}), which is now established, we thus obtain that
\begin{eqnarray}
 \Q\Big(\big|\Ccal_{\sss(1)}-Nq^*_{\sN}\big| > \omega_{\sN}\sqrt{B_{\sN}}\Big)
&\leq&
\Q\Big(\big|X_{\sN}-Nq^*_{\sN}\big|
> \omega_{\sN}\sqrt{B_{\sN}}\Big)+\Q(|\Ccal_{\sss(2)}|>\gamma_{\sN}^*)\nonumber\\
  &=& \Q(E_{\sN}^c)+\Q(|\Ccal_{\sss (2)}|> \gamma^*_{\sN})
    =o(1).
\end{eqnarray}
This completes the proof of Proposition \ref{prop-red}.
\qed

The proof of Theorem \ref{thm-cluster} follows from
Proposition \ref{prop-red}, by taking $Q=\prob_{\sN}$, and
$B_{\sN}=C \frac {\gamma_{\sN}^2 N^2}{L_{\sN}}=o(N^2)$.
For this we note that $\expec_{\sN}[X_{\sN}]=Nq_{\sN}$
follows from (\ref{ass3}). The second assumption in (\ref{ass3rep}) follows from
(\ref{ass1}) and Lemma \ref{lem-varbd} as follows:
\eq
\label{toevoeging2}
\sum_{i,j}\prob_{\sN}(i,j\,\mbox{connected})
=\sum_{i,j}\sum_{l}\prob_{\sN}(i,j\in \Ccal_{\sss (l)})
=\sum_l \expec_{\sN}[|\Ccal_{\sss (l)}|^2 I[\Ccal_{\sss (l)}\geq \gamma_{\sN}]]+o(N^2),
\en
because $\gamma_{\sN}=o(N)$. In turn:
\eq
\label{toevoeging3}
\sum_l |\Ccal_{\sss (l)}|^2 I[\Ccal_{\sss (l)}\geq \gamma_{\sN}]
=
\sum_{i,j} I[|\Ccal_{\sss (i)}|\ge \gamma_{\sN},|\Ccal_{\sss (j)}|\ge \gamma_{\sN}],
\en
so that
\eq
\label{toevoeging4}
\sum_l \expec_{\sN}[|\Ccal_{\sss (l)}|^2 I[\Ccal_{\sss (l)}\geq \gamma_{\sN}]]
=\expec_{\sN}[X^2_{\sN}]=(\expec_{\sN}[X_{\sN}])^2+\text{Var}_{\sN}(X_{\sN})
=N^2q_{\sN}^2(1+o(1)),
\en
because Lemma \ref{lem-varbd} stated that $\text{Var}_{\sN}(X_{\sN})$ is of order $N$.
\qed

\subsection{Proof of Theorem \ref{thm-clusterP}}
\label{sec-pfThmclusterP}

The proof of Theorem \ref{thm-clusterP} will be given by verifying the conditions
of Proposition \ref{prop-red} with $\mathbb{Q}=\prob$ and
$\gamma_{\sN}^*=\bar\gamma_{\sN}$ defined in
(\ref{gammabarNdef}). In order to do so
we will use results proved in
\cite{HHV05} for $\tau>3$, \cite{HHZ04a} for $\tau\in (2,3)$
and \cite{EHHZ04} for $\tau\in (1,2)$
(when we apply these results we will give more specific references).

We now turn to the proof of the theorem in question.
Because in the configuration model the nodes $1,2,\ldots,N$ are exchangeable,
    \eq
    \expec[X_{\sN}]=N\prob(|\Ccal_1|>\bar\gamma_{\sN}),
    \en
and this identifies $q_{\sN}=\prob(|\Ccal_1|>\bar\gamma_{\sN})$
(see (\ref{ass3})).
We next note that (again using that the nodes $1,2,\ldots,N$ are exchangeable),
    \eq
    \sum_{i,j} \mathbb{P}(i,j \text{ connected})
    =N(N-1)\mathbb{P}(1,2 \text{ connected})+N.
    \en
In \cite[p.~99, Equation (4.22)]{HHV05} (case $\tau>3$),
\cite[(4.96)]{HHZ04a} (for $\tau\in (2,3)$), it was shown that
    \eq
    \label{waardeq}
    \mathbb{P}(1,2 \text{ connected})
    =q^2(1+o(1)),
    \en
where  $q$ is the survival probability of the
delayed branching process $\{{\cal Z}_l\}_{l\ge 1} $.
For $\tau\in (1,2)$ we showed in \cite[Theorem 1.1]{EHHZ04}
that the graph-distance between between $1$ and $2$ is {\bf whp} either
equal to $2$ or to $3$, so in this case
(\ref{waardeq}) holds with $q=1$.

Comparing the conditions of Proposition \ref{prop-red}
and those of Theorem \ref{thm-clusterP} shows that
in order to use Proposition \ref{prop-red}, it remains to show that:
(i) $q_{\sN}=q+o(1)$, and (ii) to give a bound $B_{\sN}=o(N^2)$ on $\text{Var}(X_{\sN})$.
This is indeed so, because
$L_N\ge 2N$ follows from $\mu >2$, when $\tau>2$ or is immediate from
$\tau \in (1,2)$.
We prove (i) in Lemma \ref{lem-qNasy} and (ii) in Lemma
\ref{lem-varbdP} below.

\begin{lemma}
\label{lem-qNasy}
$q_{\sN}=q+o(1)$.
\end{lemma}

\proof We have that
    \eq
    q_{\sN}=\prob(|\Ccal_1|> \bar\gamma_{\sN})=1-\prob(|\Ccal_1|\leq \bar\gamma_{\sN}).
    \en
Using (\ref{couplingZhat}), we obtain
    \eq
    \label{fourfyftyone}
    \prob(|\Ccal_1|\leq \bar\gamma_{\sN})
    = \expec\big[\prob_{\sN}\big(\sum_{l} \hat Z_{l}^{\smallsup{1,N}}\leq \bar \gamma_{\sN}\big)\big]
    + O(\frac{\bar\gamma_{\sN}^2}{L_{\sN}}).
    \en
The coupling is described in full detail
in \cite[Section 3, p.~87]{HHV05}, whereas the bound in (\ref{couplingZhat}) and hence
(\ref{fourfyftyone})
follows from the proof of \cite[Lemma A.2.2, p.~111]{HHV05}, which holds under the
$\prob_{\sN}$-probability and is therefore true for any degree sequence.
Therefore,
    \eq
    q_{\sN}=1-\expec\big[\prob_{\sN}\big(\sum_{l} \hat Z_{l}^{\smallsup{1,N}}\leq \bar \gamma_{\sN}\big)\big]
    + O(\frac{\bar\gamma_{\sN}^2}{L_{\sN}}).
    \en
We start with $\tau \in (1,2)$.
We note that with probability 1,
    \eqalign
    \prob_{\sN}\big(\sum_{l} \hat Z_{l}^{\smallsup{1,N}}\leq \bar \gamma_{\sN}\big)
    &\leq \prob_{\sN}(\hat Z_{2}^{\smallsup{1,N}}\leq\bar\gamma_{\sN})
    \label{tau(1,2)couplingbd}
    \leq \sum_{n=1}^{\bar\gamma_{\sN}} g_n^{\smallsup{N}}=
    \sum_{i=1}^N \frac{D_i}{L_{\sN}}I[D_i\leq
    \bar\gamma_{\sN}+1]\leq \frac{(\bar\gamma_{\sN}+1)N}{L_{\sN}}.
    \enalign
Therefore, since by dominated convergence both
$\expec[\frac{\bar\gamma_{\sN}^2}{L_{\sN}}]\to 0$ and $\expec[\frac{(\bar\gamma_{\sN}+1)N}{L_{\sN}}]\to 0$, we conclude that
$q_{\sN}=1-o(1)$, when $\tau\in(1,2)$.

We next turn to $\tau\in (2,3)$ and $\tau>3$, which we treat simultaneously.
For this, we use that we can prove by coupling (see \cite[Section 3, p.~87]{HHV05}) that
    \eq
    \label{coupling}
    \prob_{\sN}(\sum_{l} \hat Z_{l}^{\smallsup{1,N}}\leq \bar\gamma_{\sN})
    =\prob_{\sN}(\sum_{l} {\cal Z}_{l}\leq \bar\gamma_{\sN})
    +O(\bar\gamma_{\sN} p_{\sN})=\prob(\sum_{l} {\cal Z}_{l}\leq\bar\gamma_{\sN})
    +O(\bar\gamma_{\sN} p_{\sN}),
    \en
where $p_{\sN}$ is the total variation distance between $\{g_n^{\smallsup{N}}\}$ and
$\{g_n\}$ given by
    \eq
    p_{\sN}=\frac 12 \sum_{n=0}^{\infty} |g_n^{\smallsup{N}}-g_n|,
    \en
and where the second equality in (\ref{coupling}) follows since the offspring distribution of
$\{{\cal Z}_{l}\}_{l\ge 0}$ does not depend on the degrees $D_1,D_2,\ldots,D_{\sN}$.

In \cite[Proposition 3.4, p.~92]{HHV05}, it is shown that for $\tau>3$, and some $\alpha_2,\beta_2>0$,
\begin{equation}
\label{proppN}
\prob(p_{\sN}>N^{-\alpha_2})\leq N^{-\beta_2}.
\end{equation}
In \cite[Remark A.1.3, p.~107]{HHV05}, the same conclusion is derived for $\tau\in (2,3)$.
Therefore,
    \eq
    q_{\sN}=1-\prob\big(\sum_{l} {\cal Z}_{l}\leq \bar \gamma_{\sN}\big)
    +O(\frac{\bar\gamma_{\sN}^2}{L_{\sN}})+O(N^{-\beta_2})+O(\bar \gamma_{\sN}N^{-\alpha_2}),
    \en
so that, in turn,
    \eq
    q_{\sN}=q-\prob\big(\bar \gamma_{\sN}< \sum_{l} {\cal Z}_{l}<\infty\big)
    +o(1).
        \en
We have that
    \eq
    \prob\big(\bar \gamma_{\sN}< \sum_{l} {\cal Z}_{l}<\infty\big)
    =(1-q)\prob\big(\sum_{l} {\cal Z}_{l}> \bar \gamma_{\sN}
    |\text{ extinction}\big).
    \en
A supercritical branching process conditioned
on extinction is a branching process with
law
    \eq
    g^*_n= (1-q)^{n-1} g_n \qquad (n\geq 1), \qquad g_0^*=1-\sum_{n\geq 1}g^*_n.
    \en
Indeed, if $F_n$ is the event that ${\cal Z}_l$ has $n$
children in the first generation, then
\begin{equation}
\prob( {\cal Z}_l\,\mbox{dies out},F_n)=
g_n
\prob(\mbox{$n$ copies of}\, {\cal Z}_l\,\mbox{die out})=(1-q)^n g_n.
\end{equation}
It is not hard to see that $g^*$ is a subcritical offspring distribution, and it clearly has finite mean.
Therefore, in particular, the total progeny has finite mean (in fact,
even exponential tails), so that by the Markov
inequality
    \eq
    \prob\big(\sum_{l} {\cal Z}_{l}> \bar \gamma_{\sN}
    |\text{ extinction}\big)\leq \bar \gamma_{\sN}^{-1} \expec[\sum_{l}
    {\cal Z}_{l}|\mbox{extinction}]
    =O(\bar \gamma_{\sN}^{-1})=o(1).
    \en
This completes the proof of Lemma \ref{lem-qNasy}.
\qed

We must also show (ii), i.e., we have to show that the variance of $X_N$
is bounded by $B_{\sN}=o(N^2)$. We will show:

\begin{lemma}
\label{lem-varbdP}
There exists $\beta>0$ such that
\eq
\label{VarbdP}
{\rm Var}(X_{\sN})=O(N^{2-\beta}).
\en
\end{lemma}

\proof We follow the proof of Lemma \ref{lem-varbd}.
We rewrite
    \eq
    \text{Var}(X_{\sN})=\text{Var}(Y_{\sN})
    =\expec\big(\text{Var}_{\sN}(Y_{\sN})\big)
    +\expec\big(\expec_{\sN}[Y_{\sN}]^2\big)-\expec[Y_{\sN}]^2.
    \en
By Lemma \ref{lem-varbd}, $\expec\big(\text{Var}_{\sN}(Y_{\sN})\big)$ is certainly bounded by
$O(N^{2-\beta})$,
and we are left to bound the second term. We start with $\tau \in (1,2)$.
We use (\ref{tau(1,2)couplingbd}) to see that {\bf whp}, and for some $\eta>0$
\begin{equation}
\prob_{\sN}(\sum_{l}\hat Z_{l}^{\smallsup{1,N}}\leq \bar\gamma_{\sN})
\leq \frac{(\bar\gamma_{\sN}+1)N}{L_{\sN}}\leq N^{-\eta}.
\end{equation}
 Therefore, using that $Y_{\sN}=\sum_{i=1}^N I[|C_i|\leq  \bar\gamma_{\sN}]$,
    \eq
   \expec\big(\expec_{\sN}[Y_{\sN}]^2\big)
    \leq\expec\big(N^2\prob^2_{\sN}(\sum_l \hat Z_{l}^{\smallsup{1,N}}\leq \bar\gamma_{\sN})\big)
    + O(N^2\frac{\bar\gamma_{\sN}^4}{L_{\sN}^2})=O( N^{2-\eta}).
    \en
Therefore, (\ref{VarbdP}) holds with $\beta=\eta$.

We next turn to $\tau\in (2,3)$ and $\tau>3$, which we treat simultaneously.
Using once more the fact that $Y_{\sN}=\sum_{i=1}^N I[|C_i|\leq  \bar\gamma_{\sN}]$
and (\ref{couplingZhat}), we get
\begin{equation}
\expec_{\sN}[Y_{\sN}]=\sum_{i=1}^N \prob_{\sN}(|C_i|\leq  \bar\gamma_{\sN})
=
\sum_{i=1}^N
 \prob_{\sN}(\sum_{l} \hat Z_{l}^{\smallsup{i,N}}\leq \bar\gamma_{\sN})+O(\frac{N\bar\gamma_{\sN}^2}{L_{\sN}}).
\end{equation}
Hence
\begin{eqnarray}
\expec\big(\expec_{\sN}[Y_{\sN}]^2\big)-\expec[Y_{\sN}]^2
&=&\sum_{i,j}
\left\{
\expec\big[
\prob_{\sN}(\sum_{l} \hat Z_{l}^{\smallsup{i,N}}\leq \bar\gamma_{\sN})
\prob_{\sN}(\sum_{l} \hat Z_{l}^{\smallsup{j,N}}\leq \bar\gamma_{\sN})
\big]
\right.\nn\\
&&\quad-
\left.
\prob(\sum_{l} \hat Z_{l}^{\smallsup{i,N}}\leq \bar\gamma_{\sN})
\prob(\sum_{l} \hat Z_{l}^{\smallsup{j,N}}\leq \bar\gamma_{\sN})
\right\}+O(\frac{N^2\bar\gamma_{\sN}^2}{L_{\sN}}).
\end{eqnarray}
Now, by (\ref{coupling}), we can replace $\prob_{\sN}(\sum_{l} \hat Z_{l}^{\smallsup{i,N}}\leq \bar\gamma_{\sN})$
by $\prob_{\sN}(\sum_{l} {\cal Z}_{l}^{\smallsup{i}}\leq \bar\gamma_{\sN})$, at the cost of an
additional error term $O({\bar\gamma}_{\sN}p_{\sN})$,
where $\{{\cal Z}_{l}^{\smallsup{i}}\}_{l\ge 1}$
for $i=1,2,\ldots,N$ are independent copies of the branching process ${\cal Z}_{l}$.
Since $\{{\cal Z}_{l}^{\smallsup{i}}\}_l$
and $\{{\cal Z}_{l}^{\smallsup{j}}\}_l$ are independent for $i\neq j$ and their law is
independent of the degree sequence, we have that
    \eq
    \prob_{\sN}(\sum_{l} {\cal Z}_{l}^{\smallsup{i}}\leq \bar\gamma_{\sN})
    \prob_{\sN}(\sum_{l} {\cal Z}_{l}^{\smallsup{j}}\leq \bar\gamma_{\sN})
    =\prob^2(\sum_{l} {\cal Z}_{l}\leq\bar\gamma_{\sN}),
    \en
so that we obtain
\begin{eqnarray}
\expec\big(\expec_{\sN}[Y_{\sN}]^2\big)-\expec[Y_{\sN}]^2
&=&
O(N^2{\bar\gamma}_{\sN}\expec[p_{\sN}])+O(\frac{N^2\bar\gamma_{\sN}^4}{L_{\sN}})=O(N^{2-\beta}),
\end{eqnarray}
by bounding the sum over $i$ by $N$,
and using (\ref{proppN}), which implies $\expec[p_{\sN}]\le N^{-(\alpha_2\wedge \beta_2)}$,
so that  choosing $0<\beta<\alpha_2\wedge \beta_2$, kills the
additional factor $\log{N}$ originating from $\bar \gamma_{\sN}$.
\qed

This concludes the proof of Theorem \ref{thm-clusterP}. We even obtain an improvement,
since $\sqrt{B_{\sN}}\leq N^{1-\beta'}$ for any $\beta'<\beta$, so that {\bf whp}
the largest cluster is in between $Nq_{\sN} \pm \omega_{\sN} N^{1-\beta'}$.
\qed

\section{Further bounds on connected components and diameter}
\label{sec-lb}
\subsection{On connected components}
\label{sec-cc}
The proof of Theorem~\ref{dnz9Th2}
is based on the following lemma. Recall that $f_k=\prob(D=k),\, k\ge 1$.
\begin{lemma}
\label{dnz9L1}
Assume the conditions of Theorem~\ref{dnz9Th2}. Suppose further that for some
$k=k_{\sN}=O(\log{N})$, and some $0<\delta<1/6$,
\whp
    \begin{equation}
    \label{dnz9_15}
    Nf_k \left(\frac{f_1(1-\delta)}{\mu_{\sN}}\right)^k \rightarrow \infty.
    \end{equation}
Then \whp the random graph contains a
connected component with $k+1$ nodes.
\end{lemma}

\noindent
{\bf Proof.}
Take $k$ such that $f_k>0$ and
consider the star-like connected component,
with one node of degree $k$ at the center and
$k$ nodes of degree $1$ at the ends
(see Figure~\ref{dnz9_f1}).
\fig{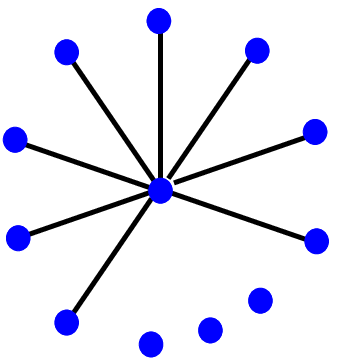}{A star-like connected component with $k+1$ nodes.}{dnz9_f1}

We will show that if the condition of the lemma
holds, then, the random graph contains
the above connected component {\bf whp}.

The main idea behind the proof is the following.
There are \whp at least $f_1(1-\delta)N$ nodes of degree one.
Hence, the probability that we connect a node of degree $k$
to $k$ nodes of degree $1$ is at least
    \begin{equation}
    \label{dnz9_2}
    \left(\frac{f_1 (1-\delta)N}{L_{\sN}}\right)^k=\left(\frac{f_1(1-\delta)}{\mu_{\sN}}\right)^k.
    \end{equation}
Since, \whp there are about $Nf_k$ nodes of degree $k$,
we have about $Nf_k$ trials to make such a $k$-star component.
The mean number of successful trials is then about
    $$
    Nf_k \left(\frac{f_1(1-\delta)}{\mu_{\sN}}\right)^k\rightarrow \infty,
    $$
by condition~(\ref{dnz9_15}).
Hence, \whp we expect to make at least one successful trial,
and find a $k$-star component.

We will now give the details of the proof.
First we define a procedure which determines the
existence of $k$-stars in the random graph.
Consider the process of pairing stubs in the graph.
We are free to choose the order in which we pair them.
Consider $D_{j_1},\dots,D_{j_{\ell_{\sN}}}$, where we abbreviate
$\ell_{\sN}=N(k)$ to be equal to the number of the nodes with degree $k$, which we call
{\it $k$-nodes} for brevity. We pair the stubs in the following order:

Let $S_{(\sss 1)}=S_{\sss 1,1},\dots,S_{\sss 1,k}$ be the stubs of node $j^*_1=j_1$.
We first pair $S_{\sss 1,1}$. If it is paired with a stub of a node of degree $1$,
then we call this {\em pairing successful} and consider the pairing of $S_{\sss 1,2}$.
If $S_{\sss 1,2}$ is paired with a stub of a node of degree $1$,
then we call the second pairing successful and consider the pairing of $S_{\sss 1,3}$,
    and so on until the first moment when one of the two following things happens.
The first case is that all stubs of node $j_1$ are paired with nodes
of degree $1$. Then we observe a $k-star$ component, we call the first
{\em trial successful} and stop.
The second case is that we come to $l<k$ such that the $l^{\rm th}$ pairing is
unsuccessful, i.e., $S_{\sss 1,l}$ is not paired with a node of degree $1$.
Then we call the whole trial unsuccessful and stop with pairing
of the stubs in $S_{\sss (1)}$.
In the later case it is possible that
$S_{\sss 1,l}$ is paired with another node in $j_2,\dots j_{\ell_{\sN}}$.
Such node can not turn into a $k$-star anymore, we call this
node, as well as node $j_1$, {\em used} and discard them in the procedure.

We define our successive trials inductively.
For any $m\ge 2$, let $j^*_m$ be the first unused node in the
sequence $j_1,\dots j_{\ell_{\sN}}$.
Then, for $j^*_m$ we use the same procedure as with $j^*_1$
to determine whether the $m^{\rm th}$ trial is successful or not.
If the trial is not successful, then node $j^*_m$ becomes used,
and if the corresponding unsuccessful pairing involves another
unused $k$-node, then we also call this node used. We always discard
all used nodes from the procedure.

We repeat these trials until we find a successful $k$-node or
until all $k$-nodes are used. The $\prob_{\sN}$-probability
that the $j^{\rm th}$ trial is successful is
    \begin{equation}
    \label{dnz9_12}
    \prod\limits_{s=0}^{k-1}
    \frac{N(1)-L_{\sN}(1,j)-s}{L_{\sN}-2jk-2s-1}\indic[N(k)\ge j]
    \indic[N(1)-L_{\sN}(1,j)\ge k],
    \end{equation}
where $N(1)-L_{\sN}(1,j)$ is the remaining number of free stubs
of the degree 1 nodes up to the moment of the $j^{\rm th}$ trial.
Let $\tau_{\sN}(k)$ be the number of trials.
Since at every unsuccessful trial the number of used nodes of degree $k$ increases
by at most two, we have $\tau_{\sN}(k)\ge \lfloor N(k)/2\rfloor$.
Instead of using all these trials, we will only use
$\lfloor \delta^2 N(k)\rfloor$ of them.
Then, after $\lfloor \delta^2 N(k)\rfloor$ trials,
there are at least $N(1)-k \lfloor \delta^2 N(k)\rfloor$
remaining free stubs attached to nodes of degree 1. Hence,
for $j\le \lfloor \delta^2 N(k)\rfloor$, (\ref{dnz9_12}) is at least
    \begin{equation}
    \label{dnz9_13}
    \left(
    \frac{N(1)-\delta^2 k N(k)}{L_{\sN}}\right)^k
    \indic[N(k)\geq j]
    \indic[N(1)-\delta^2 k N(k)\ge 0],
    \end{equation}
where we use that if $N(1)-\delta^2 k N(k)\ge 0$, then, for all
$j\leq \lfloor \delta^2 N(k)\rfloor-1$, we have that $N(1)-L_{\sN}(1,j)\ge k$.
Then, the $\prob_{\sN}$-probability that all
$\lfloor \delta^2 N(k)\rfloor$ trials are unsuccessful is at most
    \begin{equation}
    \label{dnz9_3}
    \prod\limits_{i=0}^{\lfloor \delta^2 N(k)\rfloor-1}
    \left(1-\left(
    \frac{N(1)-\delta^2 k N(k)}{L_{\sN}}\right)^k\right)\indic[N(1)-\delta^2 k N(k)\ge 0].
    \end{equation}
If we show that \whp
    \begin{equation}
    \label{dnz9_4}
    \begin{array}{rl}
    (i)&\qquad N(1)-\delta^2 k N(k)\ge (1-\delta)Nf_1,\\
    (ii)&\qquad N(k)\ge Nf_k/2,
    \end{array}
    \end{equation}
then, again {\bf whp}, the probability that all $\lfloor \delta^2 N(k)\rfloor$  trials are unsuccessful is at most
\eq
\prod_{i=0}^{N\delta^2f_k/2} \left(
1-\left(\frac{(1-\delta)Nf_1}{L_{\sN}}\right)^k
\right)  \le
\exp\left\{-(N\delta^2f_k/2+1)
\left(\frac{(1-\delta)f_1}{(L_{\sN}/N)}\right)^k
\right\}
  =o(1),
\en
due to~(\ref{dnz9_15}), where we further used that $1-x\le e^{-x}$ for $x>0$.
Hence, we are done if we prove (\ref{dnz9_4}).

For $\tau \in (1,2)$, we have by assumption $f_k>0$, for some $k\le \gamma_2^{**}$ (see \ref{gamma**}).
It follows that $k$ is bounded and so by the law of large numbers we have
\whp\!\!,
    \begin{equation}
    \label{dnz9_14}
    (1-\delta) Nf_k \le N(k)\le (1+\delta) Nf_k.
    \end{equation}
Hence, part (\ref{dnz9_4}(ii)) is clear from the lower bound in (\ref{dnz9_14}) when
$\delta\leq 1/2$.
For part (\ref{dnz9_4}(i)), we use (\ref{dnz9_14}) together with the similar bound
that \whp
    \begin{equation}
    \label{dnz9_14b}
    (1-\delta^2) Nf_1 \le N(1)\le (1+\delta^2) Nf_1.
    \end{equation}
Then
    \eq
    N(1)-\delta^2 k N(k)
    \geq (1-\delta^2) Nf_1-\delta^2(1+\delta) Nkf_k=N\big[f_1-\delta^2(f_1+(1+\delta)kf_k)\big]
    \geq Nf_1 (1-\delta),
    \en
when $\delta$ is sufficiently small, since $k$ is fixed. This completes the
proof of (\ref{dnz9_4}) when $\tau\in (1,2)$.

We turn to the case $\tau>2$. In this case $k=k_{\sN}=O(\log N)$ and
hence the law of large numbers does not apply.
Instead, we use \cite{Jans02}, which states
that a binomial random variable $X$  satisfies
\begin{equation}
\prob(|X-\expec[X]|\geq t) \leq 2e^{-\frac{t^2}{2(\expec[X]+t/3)}},
\label{binbd}
\end{equation}
for all $t>0$. We apply the above result with
$X=N(k)$, $\expec[X]=Nf_k$ and $t=\delta Nf_k$. Then we obtain for large enough $N$,
    \eq
    \prob(|N(k)-N f_k|\geq \delta f_k N)
    \le
    2e^{-\frac{\delta^2 Nf_k}{2(1+\delta/3)}}=o(1),
    \en
uniformly in $N$ and $k=k_{\sN}=O(\log{N})$.
This yields (\ref{dnz9_14}) and hence (\ref{dnz9_4}(ii)), \whps, when $\delta\leq 1/2$.
Furthermore, for $\tau>2$ we obtain because the expectation $\sum jf_j<\infty$   that
$kf_k=o(1)$ when $k\rightarrow \infty$.
Hence, \whp for large enough $N$,
    $$
    N(1)-\delta^2 k N(k)\ge (1-\delta^2)f_1N-\delta^2(1+\delta)kf_kN/2\ge(1-\delta)f_1N,
    $$
uniformly in $N$ and in $k=k_{\sN}= O(\log{N})$, for $\delta>0$
small enough, so that~(\ref{dnz9_4}(i)) holds. This completes the proof
of (\ref{dnz9_4}) for $\tau>2$.

\qed
\vskip0.5cm

\noindent
{\bf Proof of Theorem~\ref{dnz9Th2}(i)}.
We check the conditions of Lemma~\ref{dnz9L1} for $\tau \in (1,2)$.
Firstly, since $k\le \gamma^{**}$, and $\gamma^{**}$ being constant, the condition $k=k_{\sN}=O(\log{N})$,
of Lemma \ref{dnz9L1} is trivially fulfilled.
Secondly, we rewrite the left side of~(\ref{dnz9_15}) as
    \begin{equation}
    \label{dnz9_10b}
    Nf_k\left(\frac{f_1(1-\delta)}{\mu_{\sN}}\right)^k
    =f_k\left(f_1(1-\delta)\right)^k
    e^{\log{N}-k\log(\mu_{\sN})}.
    \end{equation}
Fix $0<\delta'< \delta<1/6$, and let $\varepsilon>0$ be arbitrary.
Since $L_{\sN}=D_1+\dots+D_{\sN},$ where $D_i$ is in the domain of attraction of a stable law
(\cite[Corollary 2, XVII.5, p.~578]{Fell71}), we have
    \eq
    \prob\left(\log\mu_{\sN}\le (1+\delta')\frac{2-\tau}{\tau-1}\log{N}\right)=
    \prob\left(\mu_{\sN}\le N^{(1+\delta')\frac{2-\tau}{\tau-1}}\right)=
    \prob\left(L_{\sN}\le N^{\frac{1}{\tau-1}+\frac{\delta'(2-\tau)}{\tau-1}}\right)
    \ge1-\varepsilon,
    \en
since $(2-\tau)/(\tau-1)>0$, for $\tau \in (1,2)$.
Thus, we obtain that, with probability at least $1-\varepsilon$,
    \eq
    Nf_k\left(\frac{f_1(1-\delta)}{\mu_{\sN}}\right)^k
    =f_k\left(f_1(1-\delta)\right)^k
    e^{\log{N}-k\log\mu_{\sN}}
    \geq f_k\left(f_1(1-\delta)\right)^k N^{1-(1+\delta')(1-\delta)}
    \rightarrow \infty,
    \en
for every $k\leq \gamma^{**}_2$, where $\gamma^{**}_2$ is defined in
(\ref{gamma**}). Therefore, the conditions
of Lemma \ref{dnz9L1} are fulfilled, and Theorem~\ref{dnz9Th2}(i) follows.
\qed
\vskip0.5cm

\noindent
{\bf Proof of Theorem~\ref{dnz9Th2}(ii).}
We again use Lemma~\ref{dnz9L1} and check its conditions.
Firstly, by the condition of the theorem, $k$ clearly satisfies $k=k_{\sN}=O(\log{N})$.
Secondly, we rewrite the expression in the left side of~(\ref{dnz9_15}), using~(\ref{dnz9_10}),
and with $\delta$ replaced by $\delta'$, as
    \begin{equation}
    \label{dnz9_9}
    Nf_k\left(\frac{f_1(1-\delta')}{\mu_{\sN}}\right)^k
    =L_f(k)e^{-\tau\log{k}}e^{\log{N}+k\log\left(f_1(1-\delta')/\mu_{\sN}\right)}.
    \end{equation}
Since $\tau>2$, we have by the weak law of large numbers ({\bf w.l.l.n.}), $
\mu_{\sN}\to\mu$ in probability, as $N\to\infty,$ so that
\whp $\mu_{\sN}\leq \mu/(1-\delta')$. On this event, we then obtain the lower bound
\begin{eqnarray}
\log{N}+k\log\left(\frac{f_1(1-\delta')}{\mu_{\sN}}\right)
&\geq&
\log{N}-\gamma_1^{**}\log N\cdot\log\left(\frac{\mu}{f_1(1-\delta')^2}\right)
\\
&\geq&
\log{N}\Big(1-\frac{1-\delta}{\log(\mu/f_1)}\big[\log(\mu/f_1)
-2\log(1-\delta')\big]\Big)\geq \frac{\delta}{2}\log{N},\nn
\end{eqnarray}
when $\delta'>0$ is sufficiently small.
Substituting the above lower bound in the right side of~(\ref{dnz9_9}),
we obtain from \eqref{dnz9_10} that, for sufficiently large $N$, \whp
    \eq
    Nf_k\left(\frac{f_1(1-\delta')}{\mu_{\sN}}\right)^k
    \geq L_f(k)e^{-\tau\log{k}}
    N^{\delta/2}\to\infty,\qquad \mbox{ as }N\to\infty,
    \en
where we have used that $k=k_{\sN}=O(\log{N})$ so that $e^{-\tau\log{k}}=e^{-O(\log\log n)}$ .
We conclude that the condition~(\ref{dnz9_15}) is fulfilled
with some $\delta'>0$, and we have proved Theorem~\ref{dnz9Th2}(ii).
\qed
\vskip0.5cm

\subsection{A lower bound on the diameter}
\noindent
We now prove Theorem \ref{thm-diameterfirst} which gives a lower bound on the diameter.\\
\\
{\bf Proof of Theorem \ref{thm-diameterfirst}.}
We start by proving the claim when $f_2>0$. The idea behind the proof is simple.
Under the conditions of the theorem,
one can find, \whp, a path $\Gamma(N)$ in the random graph such that
this path consists exclusively of nodes with degree $2$
and has length at least $2\alpha\log{N}$.
This implies that the diameter $D(G)$ is at least $\alpha \log{N}$,
since the above path could be a cycle.

Below we define a procedure which proves the existence
of such a path. Consider the process of pairing stubs in the graph.
We are free to choose the order in which we pair the free stubs,
since this order is irrelevant for the distribution of the random
graph. Hence, we are allowed to start with pairing the stubs of the nodes of
degree $2$.

Let $S_{\sN}(2)=(i_1,\dots,i_{N(2)})\in \NN^{N(2)}$
be the nodes of degree $2$, where we recall that
$N(2)$ is the number of such nodes.
We will pair the stubs and at the same time define
a permutation $\Pi(N)=(i^*_1,\dots,i^*_{N(2)})$ of $S_{\sN}(2)$,
and a characteristic $\chi(N)=(\chi_1,\dots,\chi_{N(2)})$ on $\Pi(N)$,
where $\chi_j$ is either $0$ or $1$.
$\Pi(N)$ and $\chi(N)$ will be defined inductively in such a way
that for any node $i^*_k\in\Pi(N)$, $\chi_k=1$, if and only if
node $i^*_k$ is connected to node $i^*_{k+1}$. Hence, if $\chi(N)$
contains a substring of at least $2\alpha\log{N}$ ones
then the random graph contains a path $\Gamma(N)$ of length at least $2\alpha\log{N}$.

We initialize our inductive definition by $i^*_1=i_1$.
The node $i^*_1$ has two stubs,
we consider the second one and pair it to an arbitrary free stub.
If this free stub belongs to another node $j\ne i^*_1$ in $S_{\sN}(2)$
then we choose $i^*_2=j$ and $\chi_1=1$, else we choose $i^*_2=i_2$,
and $\chi_1=0$. Suppose for some $1<k\le N(2)$,
the sequences $(i^*_1,\dots,i^*_k)$ and $(\chi_1,\dots,\chi_{k-1})$
are defined. If $\chi_{k-1}=1$, then one stub of $i^*_k$
is paired to a stub of $i^*_{k-1}$, and another stub of $i^*_k$ is free,
else, if $\chi_{k-1}=0$, node $i^*_k$ has two free stubs. Thus,
node $i^*_k$ has at least one free stub. We pair this stub
to an arbitrary remaining free stub. If this second stub
belongs to node $j\in S_{\sN}(2)\setminus\{i^*_1,\dots,i^*_k\}$,
then we choose $i^*_{k+1}=j$ and $\chi_k=1$, else we choose
$i^*_{k+1}$ as the first stub in $S_{\sN}(2)\setminus\{i^*_1,\dots,i^*_k\}$,
and $\chi_k=0$. Hence, we have defined $\chi_k=1$, if and only if node $i^*_k$
is connected to node $i^*_{k+1}$.

We show that \whp
there exists a  substring of ones of length at least
$2\alpha \log N$
in the first half of $\chi_N$,
i.e., in $\chi_{\frac{1}{2}}(N)=(\chi_{i^*_1},\dots,\chi_{i^*_{\lfloor N(2)/2\rfloor}})$.
For this purpose, we couple the sequence $\chi_{\frac{1}{2}}(N)$
with a sequence $B_{\frac{1}{2}}(N)=\{\xi_k\}$,
where $\xi_k$ are i.i.d.\ Bernoulli random variables
taking value $1$ with probability $f_2/(4\mu)$,
and such that $\chi_{i^*_k}\ge \xi_k$, for all $k\in\{1,\ldots,\lfloor N(2)/2\rfloor\}$,
\whp\!\!.
Indeed, for any $1\le k\le \lfloor N(2)/2\rfloor$, the $\prob_{\sN}$-probability
that $\chi_{i^*_k}=1$ is at least
    \begin{equation}
    \label{dnz10_2}
    \frac{2N(2)-C_{\sN}(k)}{L_{\sN}-C_{\sN}(k)},
    \end{equation}
where as before $N(2)$ is the total number of nodes with degree $2$, and $C_{\sN}(k)$ is
the total number of paired stubs after $k+1$ pairings.
By definition of $C_{\sN}(k)$, for any $k\le N(2)/2$, we have
    \begin{equation}
    \label{dnz10_1}
    C_{\sN}(k)=2(k-1)+1\le N(2).
    \end{equation}
Due to the {\bf w.l.l.n.} we also have that
\whp
    \begin{equation}
    \label{dnz10_3}
    N(2)\ge f_2 N/2,\qquad L_{\sN}\le 2\mu N.
    \end{equation}
Substitution of~(\ref{dnz10_1}) and~(\ref{dnz10_3}) into~(\ref{dnz10_2})
gives us that the right side of~(\ref{dnz10_2}) is at least
    $$
    \frac{N(2)}{L_{\sN}}\ge \frac{f_2}{4\mu}.
    $$
Thus, \whp we can stochastically dominate all
coordinates of the random sequence $\chi_{\frac{1}{2}}(N)$
with an i.i.d.\ Bernoulli sequence $B_{\frac{1}{2}}(N)$ of
$Nf_2/2$ independent trials with success probability $f_2/(4\mu)$. It is well known (see \cite{ER70})
that the probability of existence of a run of
$2\alpha \log{N}$ ones converges to one whenever
$$
2\alpha \log{N}
\leq (1-\varrho) \frac{\log{(Nf_2/2)}}{|\log{(f_2/(4\mu))}|},
$$
for some $0<\varrho<1$.

We conclude that \whp the sequence $B_{\frac{1}{2}}(N)$ contains
a group (and hence a substring) of $2\alpha\log{N}$ ones. Since \whp
$\chi_{\sN}\ge B_{\frac{1}{2}}(N)$, where the ordering is
componentwise, \whp the sequence $\chi_{\sN}$ also contains the same
substring of $2\alpha\log{N}$ ones,
and hence there exists a required path consisting of
at least $2\alpha \log{N}$ nodes with degree 2.
Thus, \whp the diameter is at least $\alpha \log{N}$, and we
have proved the theorem in the case that $f_2>0$.

We now complete the proof of Theorem \ref{thm-diameterfirst} when $f_2=0$
by adapting the above argument. When $f_2=0$, and since $f_1+f_2>0$, we must have that
$f_1>0$. Let $k^*>2$ be the smallest integer such that $f_{k^*}>0$.
This $k^*$ must exist, since $f_1<1$. Denote by $N^*(2)$ the total number
of nodes of degree $k^*$ of which its first $k^*-2$ stubs are connected to
a node with degree 1. Thus, effectively, after the first $k^*-2$ stubs
have been connected to nodes with degree 1, we are left with a structure
which has 2 free stubs. These nodes will replace the $N(2)$ nodes used in the
above proof. It is not hard to see that \whp $N^*(2)\geq f_2^* N/2$
for some $f^*_2>0$. Then, the argument for $f_2>0$ can be repeated, replacing
$N(2)$ by $N^*(2)$ and $f_2$ by $f_2^*$. In more detail, for any
$1\le k\le \lfloor N^*(2)/(2k^*)\rfloor$, the $\prob_{\sN}$-probability
that $\chi_{i^*_k}=1$ is at least
    \begin{equation}
    \label{dnz10_2b}
    \frac{2N^*(2)-C_{\sN}^*(k)}{L_{\sN}-C_{\sN}^*(k)},
    \end{equation}
where $C_{\sN}^*(k)$ is
the total number of paired stubs after $k+1$ pairings of the free stubs
incident to the $N^*(2)$ nodes.
By definition of $C_{\sN}^*(k)$, for any $k\le N^*(2)/(2k^*)$, we have
    \begin{equation}
    \label{dnz10_1b}
    C_{\sN}(k)=2k^*(k-1)+1\le N^*(2).
    \end{equation}
Substitution of~(\ref{dnz10_1b}), $N^*(2)\geq f_2^* N/2$ and the bound on $L_{\sN}$
in~(\ref{dnz10_3}) into~(\ref{dnz10_2b})
gives us that the right side of~(\ref{dnz10_2b}) is at least
    $$
    \frac{N^*(2)}{L_{\sN}}\ge \frac{f_2^*}{4\mu}.
    $$
Now the proof can be completed as above. We omit further details.
\qed

\subsection{A $\log\log$ upper bound on the diameter for $\tau\in (2,3)$}
\label{sec-loglog} In this section, we investigate the diameter of the
configuration model when $f_1+f_2=0$, or equivalently $\prob(D\geq
3)=1$. We assume \eqref{infvarass} for some $\tau \in (2,3)$. We
will show that under these assumptions $C_{\sF}\log\log{N}$ is an
upper bound on the diameter of $G$ for some sufficiently large
constant $C_{\sF}$ (see Theorem \ref{thm-diametersec}).

The proof is divided into two key steps. In the first, in
Proposition \ref{prop-core}, we give a bound on the diameter
of the {\it core} of the configuration model  consisting of
all nodes with degree at least a certain power of $\log{N}$.
This argument is very close in spirit to the one in \cite{RN04},
the only difference being that we have simplified the argument
slightly.
After this, in Proposition \ref{prop-periphery}, we derive a bound
on the distance between nodes with small degree and the core.
We note that Proposition \ref{prop-core} only relies on
the assumption in \eqref{infvarass}, while Proposition \ref{prop-periphery}
only relies on the fact that $\prob(D\geq 3)=1$. We start by
investigating the core of the configuration model.

We take $\sigma>\frac{1}{3-\tau}$ and define the {\it core}
$\Core$ of the configuration model to be
    \eq
    \Core=\{i: D_i\geq (\log{N})^\sigma\},
    \en
i.e., the set of nodes with degree at least $(\log{N})^\sigma$.  Then, the diameter of the core
is bounded in the following proposition:

\begin{prop}[The diameter of the core]
\label{prop-core}
For every $\sigma>\frac{1}{3-\tau}$, the diameter of
$\Core$ is bounded above by
    \eq
    \frac{2 \log\log{N}}{|\log{(\tau-2)}|}(1+o(1)) .
    \en
\end{prop}

\proof We note that \eqref{infvarass} implies that \whp
the largest degree $D_{\sss(N)}$ satisfies
    \eq
    D_{\sss(N)} \geq u_1, \quad \mbox{where} \quad u_1=N^{\frac{1}{\tau-1}}(\log{N})^{-1},
    \en
because for $N\to \infty$,
\begin{eqnarray}
\label{verdmax}
&&\prob(D_{\sss(N)} > u_1)=1-\prob(D_{\sss(N)} \leq u_1)
=
1-[F(u_1)]^{N}
\geq
1-(1-cu_1^{1-\tau})^N\nn\\
&&\qquad= 1-\left(1-c \frac{(\log{N})^{\tau-1}} {N} \right)^N \sim
1-\exp(-c(\log{N})^{\tau-1})\to 1.
\end{eqnarray}

Define
    \eq
    \Ncal{1}=\{i: D_i\geq u_1\},
    \en
so that, \whp, $\Ncal{1}\neq \varnothing$.
For some constant $C>0$, which will be specified later, and $k\ge 2$ we define recursively
    \eq
    \label{Wk-def}
    u_{k}=C\log{N} \big(u_{k-1}\big)^{\tau-2}.
    \en
Then, we define
    \eq
    \Ncal{k}=\{i: D_i\geq u_{k}\}.
    \en
We start by identifying $u_{k}$:

\begin{lemma}[Identification of $u_{k}$]
\label{lem-powersWk}
For each $k\in \N$,
    \eq
    u_{k}=C^{a_k} (\log{N})^{b_k} N^{c_k},
    \en
with
    \eq
    c_k=\frac{(\tau-2)^{k-2}}{\tau-1},\qquad
    b_k=\frac{1}{3-\tau}-\frac{4-\tau}{3-\tau}(\tau-2)^{k-1},
    \qquad
    a_k=\frac{1-(\tau-2)^{k-1}}{3-\tau}.
    \en
\end{lemma}

\proof We will identify $a_k$, $b_k$ and $c_k$ recursively. We note that $c_1=\frac{1}{\tau-1}, b_1=-1, a_1=0$.
By \eqref{Wk-def}, we can, for $k\ge 2$, relate $a_k,b_k, c_k$ to $a_{k-1},b_{k-1}, c_{k-1}$ as follows:
    \eq
    c_k=(\tau-2)c_{k-1}, \qquad b_k=1+(\tau-2)b_{k-1},
    \qquad a_k=1+(\tau-2)a_{k-1}.
    \en
As a result, we obtain
    \eqalign
    c_k&=(\tau-2)^{k-1} c_1=\frac{(\tau-2)^{k-1}}{\tau-1},\\
    b_k&= b_1(\tau-2)^{k-1}+\sum_{i=0}^{k-2} (\tau-2)^i
    =\frac{1-(\tau-2)^{k-1}}{3-\tau}-(\tau-2)^{k-1},\\
    a_k& =\frac{1-(\tau-2)^{k-1}}{3-\tau}.
    \enalign
\qed

The key step in the proof of Proposition \ref{prop-core}
is the following lemma:

\begin{lemma}[Connectivity between $\Ncal{k-1}$ and $\Ncal{k}$]
\label{lem-connNs}
Fix $k\geq 2$, and $C>4\mu/c$ (see (\ref{outgoing degree}), and (\ref{infvarass}) respectively).
Then, the probability that there exists an
$i\in \Ncal{k}$ that is not directly connected to $\Ncal{k-1}$ is
$o(N^{-\delta})$, for some $\delta>0$ independent of $k$.
\end{lemma}

\proof We note that, by definition,
    \eq
    \sum_{i\in \Ncal{k-1}} D_i\geq u_{k-1} |\Ncal{k-1}|.
    \en
Also,
    \eq
    |\Ncal{k-1}|\sim {\rm Bin}\big(N, 1-F(u_{k-1})\big),
    \en
and we have that, by \eqref{infvarass},
    \eq
    N[1-F(u_{k-1})]\geq cN (u_{k-1})^{1-\tau},
    \en
which, by Lemma \ref{lem-powersWk}, grows as a positive power of $N$, since $c_k\leq c_2
=\frac{\tau-2}{\tau-1}<\frac{1}{\tau-1}.$
Using \eqref{binbd}, we obtain that the probability that
$|\Ncal{k-1}|$ is bounded below by
$N[1-F(u_{k-1})]/2$ is exponentially small in $N$.
As a result, we obtain that for every $k$, and \whp
    \eq
    \label{stubsincidentcore}
    \sum_{i\in \Ncal{k}} D_i\geq \frac{c}{2} N(u_{k})^{2-\tau}.
    \en

We note (see e.g., \cite[(4.34)]{HHZ04a} that for any two sets of nodes $A$, $B$, we have that
    \eq
    \label{noconnAB}
    \prob_{\sN} (A \text{ not directly connected to }B)
    \leq e^{-\frac{D_AD_B}{L_{\sN}}},
    \en
where, for any $A\subseteq \{1, \ldots, N\}$, we write
    \eq
    D_A=\sum_{i\in A} D_i.
    \en
On the event where $|\Ncal{k-1}|\geq N[1-F(u_{k-1})]/2$
and where $L_{\sN}\leq 2\mu N$, we then obtain by \eqref{noconnAB}, and Boole's inequality
that the $\prob_{\sN}$-probability that
there exists an $i\in \Ncal{k}$ such that $i$ is not directly
connected to $\Ncal{k-1}$ is bounded by
    \eq
    N e^{-\frac{u_{k} Nu_{k-1}[1-F(u_{k-1})]}{2L_{\sN}}}
    \leq N e^{-\frac{cu_{k} (u_{k-1})^{2-\tau}}{4\mu}}
    =N^{1-\frac{cC}{4\mu}},
    \en
where we have used \eqref{Wk-def}. Taking $C>4\mu/c$ proves the claim.
\qed

We now complete the proof of Proposition \ref{prop-core}.
Fix
    \eq
    k^*= \frac{2 \log\log{N}}{|\log{(\tau-2)}|}.
    \en
As a result of Lemma \ref{lem-connNs}, we have {\bf whp} that the diameter
of $\Ncal{k^*}$ is at most $2k^*$,
because the distance between any node in
$\Ncal{k^*}$ and the node with degree $D_{(\sN)}$ is at most $k^*$. Therefore, we are done when we can show that
    \eq
    \Core\subseteq \Ncal{k^*}.
    \en
For this, we note that
    \eq
    \Ncal{k^*}=\{i: D_i \geq u_{k^*}\},
    \en
so that it suffices to prove that $u_{k^*}\geq (\log{N})^{\sigma}$,
for any $\sigma>\frac1{3-\tau}$. According to Lemma \ref{lem-powersWk},
    \eq
    u_{k^*}=C^{a_{k^*}} (\log{N})^{b_{k^*}} N^{c_{k^*}}.
    \en
Because for $x\to\infty$, and $2<\tau<3$,
\eq
\label{toegerard1}
x(\tau-2)^{\frac{2\log{x}}{|\log{(\tau-2)}|}}
=x\cdot x^{-2}=o(\log{x}),
\en
we find with $x=\log{N}$ that
\eq
\label{toegerard2}
\log{N}\cdot(\tau-2)^{\frac{2\log{\log{N}}}{|\log{(\tau-2)}|}}
=o(\log{\log{N}}),
\en
implying: $N^{c_{k^*}}=(\log{N})^{o(1)}$, $(\log{N})^{b_{k^*}}= (\log{N})^{\frac{1}{3-\tau}+o(1)}$,
and $C^{a_{k^*}}=(\log{N})^{o(1)}$. Thus,
    \eq
    \label{Wkcore}
    u_{k^*}=(\log{N})^{\frac{1}{3-\tau}+o(1)},
    \en
so that, by picking $N$ sufficiently large, we can make $\frac{1}{3-\tau}+o(1)\leq \sigma$.
This completes the proof of Proposition \ref{prop-core}.
\qed
\vskip0.5cm

Define
\eq
\label{def-C(m)}
C(m,\vep)=\left(\frac{\tau-2}{3-\tau}+1+\vep\right)/\log m,
\en
where $\vep>0$ and $m\ge 2$ is an integer.

\begin{prop}[The maximal distance between the periphery and the core]
\label{prop-periphery}
Assume that $\prob(D\geq m+1)=1$ for some $m\geq 2,$ and take $\vep>0$.
Then, \whp, the maximal distance between
any node and the core is bounded from above by
$C(m,\vep)\log\log{N}$.
\end{prop}

\proof We start from a node $i$ and will show that
the probability that the distance between $i$ and $\Core$
is at least $C(m,\vep)\log\log{N}$ is $o(N^{-1})$.
This proves the claim. For this, we explore the neighborhood of
$i$ as follows. From $i$, we connect the first $m+1$ stubs (ignoring
the other ones). Then, successively, we connect the first $m$
stubs from the closest node to $i$ that we have connected to and
have not yet been explored. We call the arising process when we
have explored up to distance $k$ from the initial
node $i$ the {\it $k$-exploration tree}.

When we never connect two stubs between nodes we have connected to,
then the number of nodes we can reach in $k$ steps is {\it precisely}
equal to $(m+1)m^{k-1}$. We call an event where a stub on
the $k$-exploration tree connects to a stub incident to  a node in the  $k$-exploration tree a {\it collision}.
The number of collisions in the $k$-exploration tree is the number of
cycles or self-loops in it. When $k$ increases, the probability
of a collision increases. However, for $k$ of order $\log\log{N}$, the probability that  more than {\it two}
collisions occur in the $k$-exploration tree is small, as we will prove now:

\begin{lemma} [Not more than one collision]
\label{lem-collision}
Take $k=\lceil C(m,\vep) \log\log{N}\rceil$. Then, the \\$\prob_{\sN}$-probability
that there exists a node of which the $k$-exploration tree has at
least two collisions, before hitting the core $\Core$, is bounded by $(\log{N})^d L_{\sN}^{-2}$,
for $d=4C(m,\vep)\log{(m+1)}+2\sigma$.
\end{lemma}

\proof For any stub in the $k$-exploration tree, the probability that it will
create a collision before hitting the core is bounded above by $(m+1)m^{k-1} (\log{N})^{\sigma}L_{\sN}^{-1}$.
The probability that
two stubs will both create a collision is, by similar arguments, bounded above
by $\big[(m+1)m^{k-1}(\log{N})^{\sigma} L_{\sN}^{-1}\big]^{2}.$ The total number of possible pairs of
stubs in the $k$-exploration tree is bounded by
$$
[(m+1)(1+m+\ldots+m^{k-1})]^2\le[(m+1)m^{k}]^2,
$$
 so that
by Boole's inequality, the probability that the $k$-exploration tree has at least two
collisions is bounded by
    \eq
    \big[(m+1)m^{k}\big]^4(\log{N})^{2\sigma} L_{\sN}^{-2}.
    \en
When $k=C(m,\vep) \log\log{N}$,
we have that $\big[(m+1)m^{k}\big]^4(\log{N})^{2\sigma}\leq (\log{N})^d$,
where\\ $d=4C(m,\vep)\log{(m+1)}+2\sigma$.
\qed
\vskip0.5cm

\noindent
Lemma \ref{lem-collision} is interesting in its own right. For example, we will
now use it together with Theorem \ref{dnzTh1intro}(i)
to prove Theorem \ref{dnzTh1intro}(ii):\\
{\it Proof of Theorem \ref{dnzTh1intro}(ii).} By Lemma \ref{lem-collision},
there are at most 2 collisions in the $k$-exploration tree from any vertex $i\in \{1, \ldots, N\}$
before hitting the core. As a result,
for any $i$, we have that the $k$-exploration tree contains at least
$\min\{(m-1) m^{k}, (\log{N})^{\sigma}\}$ stubs. When $k=C(m,\vep)\log\log{N}$,
we have that $(m-1) m^{k}\gg \log{N}$, so that the $k$-exploration tree contains at least
$K \log{N}$ stubs for some large enough $K>0$. By Proposition \ref{thm-cs},
the connected component of $i$ has \whp at least $\vep L_{\sN}$ edges,
and, in turn, by Lemma \ref{lem-smnod, largestubs}, at least $\eta N$
nodes. By Theorem \ref{dnzTh1intro}(i)
(which has already been proved in Section \ref{sec-conn} and which
applies, since $\prob(D\geq 3)=1$ and $\mu\geq 3>2$ when
$\prob(D\geq 3)=1$), we have that the size of the complement
of the largest connected component is bounded \whp.
Therefore, we must have that $i$ is part of the giant component.
Since this is true for {\it every} $i\in \{1, \ldots, N\},$
we obtain that the giant component must have size $N$, so that the
random graph is connected.
\qed
\vskip0.5cm

\noindent
Finally, we show that for $k=C(m,\vep) \log\log{N}$,
the $k$-exploration tree will, \whp connect to the $\Core$:
\begin{lemma} [Connecting the exploration tree to the core]
\label{lem-connpercore}
Take $k=C(m,\vep) \log\log{N}$. Then, the probability
that there exists an $i$ such that the distance of $i$
to the core is at least $k$ is $o(N^{-1})$.
\end{lemma}

\proof Since $\mu<\infty$ we have that $L_{\sN}/N\sim \mu$. Then, by
Lemma \ref{lem-collision}, the probability that there
exists a node for which the $k$-exploration tree has at least
2 collisions before hitting the core is $o(N^{-1})$. When the $k$-exploration tree
from a node $i$ does {\it not} have two collisions, then
there are at least $(m-1) m^{k-1}$ stubs in the $k^{\rm th}$
layer that have not yet been connected. When $k=C(m,\vep) \log\log{N}$
this number is at least equal
to $(\log{N})^{C(m,\vep)\log{m}+o(1)}$.
Furthermore, the number of stubs incident to the core
$\Core$ is stochastically bounded from below by $(\log{N})^{\sigma}$ times
a binomial distribution with parameters $N$ and success probability
$\prob(D_1\ge (\log{N})^{\sigma})$. The expected number of stubs incident to $\Core$ is
therefore at least $N(\log{N})^{\sigma}\prob(D_1\ge (\log{N})^{\sigma})$ so that
{\bf whp} the number of stubs incident to
$\Core$ is at least (by \eqref{binbd})
\eq
\label{laatste toevoeging}
\frac12 N(\log{N})^{\sigma}\prob(D_1\ge (\log{N})^{\sigma})\ge \frac{c}{2}
N(\log{N})^{\frac{2-\tau}{3-\tau}}.
\en
By \eqref{noconnAB}, the probability that we
connect none of the stubs of the $k$-exploration tree
to one of the stubs incident to $\Core$ is bounded by
    \eq
    \exp\left\{ - \frac{cN(\log{N})^{\frac{2-\tau}{3-\tau}+C(m,\vep)\log{m}} }{2L_{\sN}}\right\}\leq
    \exp\left\{-\frac{c}{4\mu}(\log{N})^{\frac{2-\tau}{3-\tau}+C(m,\vep)\log{m}}\right\}
    =o(N^{-1}),
    \en
    because \whp $L_{\sN}/N\le 2\mu$, and since $\frac{2-\tau}{3-\tau}+C(m,\vep)\log{m}=1+\vep$.
\qed
\vskip0.5cm

Propositions \ref{prop-core} and \ref{prop-periphery} prove that \whp
the diameter of the configuration model is bounded above by
$C_{\sF} \log\log{N}$, where
    \eq
    C_{\sF}= \frac{2}{|\log{(\tau-2)}|}+\frac{2(\frac{\tau-2}{3-\tau}+1+\vep)}{\log{m}}.
    \en
This completes the proof of Theorem \ref{thm-diametersec}.
\qed

\subsection*{Acknowledgement}
The work of RvdH and DZ was supported in part by
Netherlands Organisation for Scientific Research (NWO).
This work was performed in part at
the Mittag-Leffler Institute during a visit of all three authors
in the fall of 2004.


\begin{thebibliography}{99}
\bibitem{ACL01a}
W.~Aiello, F.~Chung and L.~Lu.
\newblock A random graph model for power law graphs.
\newblock {\em Experiment. Math.,} {\bf 10}(1), 53--66, (2001).


\bibitem{AB02}
R.~Albert and A.-L.~Barab\'asi.
\newblock Statistical mechanics of complex networks.
\newblock {\em Rev.\ Mod.\ Phys.} {\bf 74}, 47-97, (2002).


\bibitem{AL04}
R.A.~Arratia and T.M.~Liggett.
\newblock  How likely is an i.i.d.\ degree sequence to be graphical?
\newblock {\em Annals of Appl. Probab.,} {\bf 15}, 652-670, (2005).

\bibitem{BA99}
A.-L.~Barab\'asi and R.~Albert.
\newblock Emergence of scaling in random networks.
\newblock {\em Science} {\bf 286}, 509--512, (1999).

\bibitem{Bara02}
A.-L. Barab\'asi.
\newblock {\em Linked, The New Science of Networks}.
\newblock Perseus Publishing, Cambridge, Massachusetts, (2002).

\bibitem{Boll01}
B. Bollob\'as.
\newblock \emph{Random Graphs,} 2nd edition,
\newblock Academic Press, (2001).

\bibitem{BDM-L05}
T.~Britton, M.~Deijfen, and A.~Martin-L\"of.
\newblock Generating simple random graphs with arbitrary degree sequences,
\newblock preprint, www.math.su.se/~mia (2005).



\bibitem{CL02a}
F.~Chung and L.~Lu.
\newblock The average distances in random graphs with given expected degrees.
\newblock PNAS, {\bf 99}(25), 15879--15882, (2002).

\bibitem{CL02b}
F.~Chung and L.~Lu.
\newblock Connected components in random graphs with given expected degree sequences,
\newblock {\em Annals of Combinatorics}, {\bf 6}, 125--145, (2002).


\bibitem{CL04}
F.~Chung and L.~Lu.
\newblock The small world phenomenon in hybrid power law graphs,
\newblock In {\em Complex networks,}
{\em Lecture notes in Physics,} {\bf 650}, 89-104, Springer, Berlin, (2004).


\bibitem{CH03}
R. Cohen and S. Havlin.
\newblock Scale free networks are ultrasmall,
Physical Review Letters 90, 058701 (2003).

\bibitem{DGM02}
S.N.~Dorogovtsev, A.V.~Goltsev and J.F.F.~Mendes.
\newblock Pseudofractal scale-free web.,
\newblock {\em Phys. Rev. E} {\bf 65}, 066122, (2002).

\bibitem{DMS03}
S.N.~Dorogovtsev, J.F.F.~Mendes and A.N. Samuhkin.
\newblock Metric structure of random networks,
\newblock {\em Nucl. Phys. B.} 653, 307, (2003).

\bibitem{ER70}
P.~Erd\"os and A.~R\'enyi.
\newblock  On a new law of large numbers,
\newblock  {\em J.\ Analyse Math.} {\bf 23}, 103--111, (1970).

\bibitem{FFF99}
C.~Faloutsos, P.~Faloutsos and M.~Faloutsos.
\newblock On power-law relationships of the internet topology,
\newblock {\em Computer Communications Rev.}, {\bf 29}, 251-262, (1999).

\bibitem{EHH06}
H. van den Esker, R. van der Hofstad and G. Hooghiemstra.
\newblock Universality for the distance in finite variance random graphs,
\newblock preprint (2006).

\bibitem{EHHZ04}
H. van den Esker, R. van der Hofstad, G. Hooghiemstra and D.~Znamenski.
\newblock Distances in random graphs with infinite mean degrees,
\newblock Extremes {\bf 8}, 111-141, 2006.

\bibitem{Fell71}
W.~Feller.
\newblock {\em An Introduction to Probability Theory and Its Applications.}
Volume II, 2nd edition,
\newblock John Wiley and Sons, New York, (1971).


\bibitem{HHV05}
R. van der Hofstad, G. Hooghiemstra and P. Van Mieghem.
\newblock Distances in random graphs with finite variance degrees.
\newblock {\em Random Structures and Algorithms} {\bf 26}, 76-123, 2005.

\bibitem{HHZ04a}
R. van der Hofstad, G. Hooghiemstra and D.~Znamenski.
\newblock Distances in random graphs with finite mean and infinite variance degrees,
\newblock preprint (2006).


\bibitem{Jans02}
S.~Janson.
\newblock On concentration of probability,
\newblock Contemporary Combinatorics, ed. B. Bollob\'as,
Bolyai Soc.\ Math.\ Stud.\ {\bf 10}, J\'anos Bolyai Mathematical Society, Budapest, 289-301, (2002).


\bibitem{MR95}
M.~Molloy and B.~Reed.
\newblock A critical point for random graphs with a given degree sequence,
\newblock {\em Random Structures and Algorithms}, {\bf 6}, 161-179, (1995).


\bibitem{MR98}
M.~Molloy and B.~Reed.
\newblock The size of the giant component of a random graph with a given degree sequence,
\newblock {\em Combin.\  Probab.\ Comput.}, {\bf 7}, 295-305, (1998).


\bibitem{Newm03}
M.E.J. Newman.
\newblock The structure and function of complex networks,
\newblock {\em SIAM Rev.} {\bf 45}(2), 167--256, (2003).


\bibitem{RN04}
H.~Reittu and I.~Norros.
\newblock On the power law random graph model of massive data networks,
\newblock {\em Performance Evalution}, {\bf 55} (1-2), 3-23, (2004).


\bibitem{Stro01}
S.~H. Strogatz.
\newblock Exploring complex networks,
\newblock {\em Nature}, {\bf 410}(8), 268--276, March (2001).


\bibitem{Watt99}
D.~J. Watts.
\newblock {\em Small Worlds, The Dynamics of Networks between Order and
Randomness},
\newblock Princeton University Press, Princeton, New Jersey, (1999).


\end{thebibliography}
\end{document}